\documentclass[a4paper,11pt]{article}
\usepackage{amsmath,amssymb}
\usepackage{datetime}
\usepackage{color}
\usepackage{verbatim}
\usepackage{amsthm}
\usepackage{amssymb}
\usepackage{amsmath}
\usepackage{graphicx}
\usepackage{appendix}
\usepackage{color} 
\usepackage{jheppub}
\usepackage{amssymb}
\usepackage{amsmath}
\usepackage{mathdots}
\usepackage{verbatim}
\usepackage[utf8]{inputenc}
\usepackage[usenames,dvipsnames,svgnames,table]{xcolor}
\usepackage[all]{xy}
\usepackage{graphicx}
\usepackage{mathrsfs}
\usepackage{hyperref}

\newtheorem{thm}{Theorem}[section]
\newtheorem{lem}[thm]{Lemma}
\newtheorem{prop}[thm]{Proposition}
\newtheorem{cor}[thm]{Corollary}

\newtheorem{defn}{Definition}[section]

\newtheorem{re}{Remark}[section]

\newcommand\TM{Teichm\"uller~}

\newcommand{\ub}{\mathrm{b}}

\newcommand{\vm}{{v^{(-)}}}
\newcommand{\vp}{{v^{(+)}}}

\newcommand{\hk}{{\hat k}}

\newcommand{\ka}{k_\alpha}
\newcommand{\kb}{k_\beta}
\newcommand{\hka}{{\hat k}_\alpha}
\newcommand{\hkb}{{\hat k}_\beta}
\newcommand{\la}{l_\alpha}
\newcommand{\lb}{l_\beta}
\newcommand{\tka}{{\tilde k}_\alpha}
\newcommand{\tkb}{{\tilde k}_\beta}
\newcommand{\tla}{{\tilde l}_\alpha}
\newcommand{\tlb}{{\tilde l}_\beta}

\newcommand{\Hom}{\mathrm{Hom}}

\newcommand\bC{\mathbb{C}}
\newcommand\cA{\mathcal{A}}
\newcommand\cB{\mathcal{B}}
\newcommand\cT{\mathcal{T}}

\def\to{\longrightarrow}

\title{Heisenberg and Drinfeld doubles of $U_q(gl(1|1))$ and $U_q(osp(1|2))$ super-algebras}

\author[a,b]{Nezhla Aghaei}
\author[c]{M. K. Pawelkiewicz} 
\affiliation[a] {Section of Mathematics, University of Geneva, Rue du Conseil-Général 7-9, 1205 Geneva, Switzerland, Switzerland,}
\affiliation[b] {Centre for Quantum Mathematics
	Department of Mathematics and Computer Science (IMADA), Southern Denmark University (SDU), Campusvej 55, 5230 Odense, Denmark.}
\affiliation[c]{ CentraleSupélec, LGPM, 3 Rue Joliot Curie, Plateau de Moulon, 91190 Gif-sur-Yvette, France.}

\emailAdd{
	Nezhla.Aghaei@unige.ch
	~michal.pawelkiewicz@centralesupelec.fr}

\abstract{
	We study the Heisenberg double and the Drinfeld double of the Borel half of $U_q(gl(1|1))$ and of the $U_q(gl(1|1))$ when $q$ is a root of unity. We also study the Borel half of $U_q(osp(1|2))$ for both cases when $q$ is a root of unity and when it is not.
	We prove the isomorphism between the Heisenberg doubles and the handle algebras, which is missing in the literature, and extend the isomorphism to the graded Heisenberg doubles and the handle algebras in the context of the $\mathbb{Z}_2$-graded generalisation of Alekseev-Schomerus combinatorial quantisation of Chern-Simons theory \cite{AGPS1,AGPS2}, as well as illustrate it on the example of the Heisenberg double of the $U_q(gl(1|1))$ Hopf algebra for $q$ being a root of unity. In addition, we generalise an isomorphism between the Drinfeld double and the loop algebra from the Alekseev-Schomerus combinatorial quantisation to the graded setting.
}

\begin{document} 

\maketitle
\flushbottom

\section{Introduction}\label{chapter1}

\textit{Hopf algebra} $U_q(g)$ have been initially defined by Drinfeld \cite{Drinfeld:1985rx} and Jimbo \cite{Jimbo:1985zk} for finite-dimensional simple Lie algebras. Those $U_q(g)$ constitute a $1$-parameter deformation of the universal enveloping algebras $U(g)$ of Lie algebras $g$, and this deformation extends to the representation theory --- the irreducible finite-dimensional representations of universal enveloping algebras deform to irreducible representations of Hopf algebras when the deformation parameter $q$ is not a root of unity.

The methods coming from the representation theory of Hopf algebra have found a wide range of applications to the mathematical and theoretical physics (\cite{review} and its references). Readers who are not familiar with Hopf algebras we refer to \cite{Majidbook,Sw}. 
Hopf algebras are relevant in conformal field theory, where the algebras of screening charges and vertex operators satisfy the relations of q-deformed Lie algebras \cite{Witten:1989rw}. Moreover, in the context of quantum integrable systems a method of obtaining scattering matrices \cite{zam} has been developed which employs Hopf algebras. These systems admit an $R$-matrix which satisfies the \textit{Yang-Baxter equation}~\cite{Baxter:1972hz, Yang:1967bm}
\begin{align}
R_{12}R_{13}R_{23} = R_{23}R_{13}R_{12}.
\end{align}

There exists a systematic procedure to obtain solutions to the Yang-Baxter equation, which is based on the mathematical construction known under the name of \textit{Drinfeld double construction} \cite{DD,Drinfel'd1988,Majid,Majid0,Majidbook} and sometimes called \textit{quantum double} in the literature. It allows to associate a new, quasi-triangular Hopf algebra $D(\mathcal{A})$, i.e. a Hopf algebra that admits an universal $R$-matrix, which satisfies the Yang-Baxter equation, to an arbitrary given Hopf algebra $\mathcal{A}$. The Drinfeld double construction has been generalised to $\mathbb{Z}_2$-graded Hopf algebras in \cite{Gould}, and extended within the category theory framework to symmetric tensor categories. The generalisation to the setting of strictly braided categories has been obtained using double-bosonisation construction \cite{Majid:1999,Majid:1990,Majid:1997, aziz}.

Another existing double construction $H(\mathcal{A})$, called the \textit{Heisenberg double construction} \cite{lu,Ka3}, admits a canonical element $W$ that satisfies not the Yang-Baxter equation, but rather a \textit{pentagon equation}~\cite{Maillet}
\begin{align}
W_{12}W_{13}W_{23}=W_{23}W_{12}.
\end{align}
The Heisenberg double construction has been generalised in categorical theory setting to Hopf algebroid instead of Hopf algebras \cite{Xu:1999eh}.

Using Heisenberg doubles one can construct Drinfeld doubles (and as a consequence - also design their representations), as one can embed the elements of the Drinfeld double into a tensor square of Heisenberg doubles \cite{lu, Ka3}. Using this algebra map one can construct the universal R-matrix in terms of the canonical elements of the Heisenberg doubles.

\subsection*{Application of Heisenberg doubles and Drinfeld doubles in physics}
The Heisenberg and Drinfeld doubles for non-graded Hopf algebra find many applications in the non-supersymmetic physics setting, including Liouville theory, 3-dimensional quantum gravity, loop gravity and integrable systems. A limited number of those has been generalised to the $\mathbb{Z}_2$-graded setting.

\begin{itemize}
 
\item  Mathematical structure behind 3-dimensional quantum gravity has attracted strong interest in mathematics --- in particular, the relation of 3d gravity to the theory of link and knot invariants and the role of quantum groups and ribbon categories in the quantisation of it --- see \cite{Carlip:2004ba} for an overview. As Chern-Simons gauge theory is deeply connected to 3d quantum gravity, the former allowed one to employ gauge theoretical concepts and methods in the quantisation of the latter. Quantum group symmetries play an important role in the quantisation of gravity, as they appear in the Hamiltonian quantisation formalisms, e.g. in the combinatorial approach to the quantisation of Chern-Simons theory due to Alekseev-Schomerus \cite{Alekseev:1994pa,Alekseev:1994au,Alekseev:1995rn, Alekseev:1996ns,VS} and Buffenoir and Roche \cite{BR95, BR96}.
 These structures are well-understood in the quantum Chern-Simons theories with compact gauge groups. 
 
One started the combinatorial quantisation from the (classical) Fock and Rosly's description \cite{Fock:1998nu} of the Poisson structure on the moduli spaces of flat $G$-connections on Riemann surfaces. Imposing the flatness requirement requires applying Dirac's constraint quantisation programme to Poisson structure on the moduli spaces of flat connections. The Poisson–Lie symmetries constitute the gauge symmetries of the classical theory and the classical $r$-matrix arising from Drinfeld doubles automatically satisfies the constraints to provide a phase space for (2+1)-dimensional gravity. As a consequence, the quantum double is believed to play an important role in (2+1)-dimensional quantum gravity.

The Alekseev-Schomerus combinatorial quantisation of Chern-Simons theory proceeds in two steps. The first is the construction of the so-called graph algebra and its irreducible representations. The graph algebra is the quantum counterpart of Fock and Rosly’s Poisson algebra associated with a  first homotopy group on a Riemann surface $\Sigma$. Because of the topological nature of the theory, observables associated to any paths can be obtained from the ones corresponding to the generators of the first homotopy group of the surface.

In particular, the graph algebra for a Riemann surface of genus $g = 1$ and number of punctures $n = 0$ is known under the name of the \textit{handle algebra} and is isomorphic to the Heisenberg double of the Hopf algebra associated to the gauge group $G$. Imposition of the gauge group symmetries implemented by a Hopf algebra $U_q(g)$ then allows one to obtain the observable algebras and the Hilbert spaces of physical states.

Moreover, of interest is the graph algebra for a Riemann surface of genus $g = 0$ and number of punctures $n = 1$, i.e. a 1-punctured sphere $\Sigma_{0,1}$, known as the \textit{loop algebra} (or the monodromy algebra). Those loop algebras appear in the context of Chern-Simons theory \cite{Alekseev:1994pa, Alekseev:1994au, Alekseev:1995rn} as well as the current algebras on periodic lattice \cite{Alekseev:1996jz}. The semi-direct product of loop algebra with the Hopf algebra $U_q(g)$, called a \textit{gauged loop algebra}, incorporates the action of the gauge group on the algebra of operators of the theory and a certain commutant within it corresponds to the algebra of observables, i.e. physically measurable operators. In the non-graded setting it has been proven that the gauged loop algebra is isomorphic to the Drinfeld double of the Hopf algebra $U_q(g)$. 

The combinatorial quantisation has been extended to $\mathbb{Z}_2$-graded, (possibly) non-semisimple quasi-triangular Hopf algebras in \cite{AGPS1, AGPS2}. In addition, it has been shown that the combinatorial quantisation of Chern–Simons theory for the Drinfeld double $D(\mathcal{A})$ is equivalent to Kitaev's lattice model \cite{Kitaev:1997wr} for a finite-dimensional semisimple Hopf algebra $\mathcal{A}$. The equivalence of Kitaev models to combinatorial models shows that the former can indeed be interpreted as a Hopf algebra valued lattice gauge theory \cite{Meusburger,MeusburgerandWise}.

\item Another application of Drinfeld doubles is in 2-dimensional physics, in particular in integrable systems. The deeper understanding of integrable systems afforded by quantum groups has allowed the construction of new integrable models including spin chains (for more details c.f. \cite{book}). 
The quantum double construction can be used to obtain a quasi-triangular Hopf algebra associated with any classical Lie algebra and its universal $R$-matrix. The universal $R$-matrix together with the irreducible representations of the quantum double provide a family of solution to Yang-Baxter equation (possibly with continuous spectral parameters). When $q$ is a root of unity, those irreducible representation are associated with spin chain Hamiltonians.

\item Heisenberg doubles and Drinfeld doubles appear not only in the context of finite-dimensional Hopf algebras, but also when the Hopf algebras are countably infinite or uncountably infinite. Physically, one can understand it as the associated physical system having, instead of a finite number of states, either countably infinite or continuous spectrum. While countably infinite spectrum is present in integrable systems, the continuous spectrum appears in the context of conformal field theory (CFT). In physics the Liouville equation defines one parameter family of models in CFT, which usually is identified with a model of 2-dimensional quantum gravity. 
The space of conformal blocks of the Liouville theory with its mapping class group representation is isomorphic to the space of states obtained by quantising the \TM spaces of Riemann surfaces due to Verlinde and proved by Teschner \cite{Teschner:2002vx,Teschner:2003em,Teschner:2005bz, teschner, T03}. One aspect of the connection of the Liouville equation to the mapping class dynamics of the Teichmüller theory has been considered in the context of a specific discretised version of the Liouville equation both on classical and quantum levels \cite{Faddeev:2000if, Faddeev:2008xy, Ka0}.

Quantum Teichmüller theory was studied in \cite{Ka1,Ka2, Ka4,Chekhov:1999tn}. Heisenberg double appears in quantum Teichmüller theory, where the canonical element $W$ satisfying the pentagon relation plays the role in one of the elementary moves as a mapping class group generator of \TM theory --- the so-called \textit{flip operator} .

\item Heisenberg and Drinfeld doubles appeared also in the context of quantum computing, which leverages quantum mechanical bits (qbits) to perform computation. Yang-Baxter equation has been applied to the problem of quantum circuit compression \cite{Peng:2022} (i.e. in the problem of reducing the number of quantum gates in a quantum circuit) and that it solutions is a universal gate for quantum computing \cite{Kauffman:2004}.

The pentagon equation has also been used in the context of topological quantum computing \cite{Kitaev:2006, Kasirajan:2021} as well as in the context of quantum circuit compression \cite{Aravanis:2022shq}, where the element $W$ sasisfing the pentagon equation can express circuits with non-local interactions (i.e. those for non-adjacent qubits) can be expressed as circuits with local ones.

\end{itemize}

\subsection*{Outline of the paper}

The paper begins with a self-contained summary of $\mathbb{Z}_2$-graded Heisenberg and Drinfeld doubles, with both basis-independent and basis-dependent, i.e. more convenient for physicists, descriptions. In the Section \ref{chapter2} we first re-call the definition of Heisenberg double of a $\mathbb{Z}_2$-graded Hopf algebra $\cA$ as a smash product algebra stemming from a left Hopf algebra action. Afterwards we provide a definition equivalent to the one above in a basis dependent way --- in order to do that, we choose a basis of graded Hopf algebra $\cA$ and a basis of a dual algebra $\cA^*$ which are dual to each other. We define the canonical element $W$ which satisfies the graded pentagon relation. Finally, we define the handle algebra appearing in the super generalisation \cite{AGPS1, AGPS2} of the combinatorial quantisation of Chern-Simons theories by Alekseev-Schomerus and prove that the graded Heisenberg double is isomorphic to the handle algebra and provide the exhaustive proof of this statement. This is a novel result, as the proof for non-graded algebras was also missing in the past literature concerning the non-graded combinatorial quantisation.

In the Section \ref{chapter3} we re-call the definition of the Drinfeld double $D(\cA)$ of a $\mathbb{Z}_2$-graded Hopf algebra $\cA$ and its associated $R$-matrix. As the Heisenberg double is not a Hopf algebra --- in opposition to the Drinfeld double, which is a quasi-triangular Hopf algebra --- they can be related to each other at most on the level of algebras. In fact, the non-graded relation that the Drinfeld double can be embedded into a tensor square of Heisenberg algebras extends to the graded case. Additionally, we define the loop algebra and its semi-direct product appearing in the super generalisation \cite{AGPS1, AGPS2} of the combinatorial quantisation of Chern-Simons theories by Alekseev-Schomerus \cite{Alekseev:1995rn} and show that the graded Drinfeld double is isomorphic to it, generalising the proof by Nill \cite{Nill:1996dv}.

The following sections constitute the main part of the novel content. In the Section \ref{sec:borel-half-gl11-heisenberg-drinfeld} we construct finite-dimensional Heisenberg and Drinfeld doubles of the Borel half algebra of $U_q(gl(1|1))$, when $q$ is a root of unity. In the Section \ref{sec:gl11-heisenberg} we construct the finite-dimensional Heisenberg double of $U_q(gl(1|1))$ at root of unity and discuss its connection to the handle algebra of the $GL(1|1)$ Chern-Simons theory. In the Section \ref{sec:borel-half-osp-heisenberg-drinfeld} we construct finite-dimensional Heisenberg and Drinfeld doubles of the Borel half algebra of $U_q(osp(1|2))$, when $q$ is a root of unity, and in the Section \ref{sec:borel-half-osp-discrete-heisenberg-drinfeld} --- Heisenberg and Drinfeld doubles of the Borel half algebra of $U_q(osp(1|2))$, when $q$ is not a root of unity.

Finally, in the outlook \ref{outlook}, we motivate the interest in the Heisenberg and Drinfeld constructions for the $\mathbb{Z}_2$-graded Hopf algebras spanned by nondenumerably many basis elements. In addition, in the Appendix \ref{appendix-sl2} we provide the construction of Heisenberg and Drinfeld doubles applied to the Borel half algebra of $U_q(sl(2))$, when $q$ is not a root of unity (which is a standard, well known application) for pedagogical reasons and illustration purposes as well as to fix conventions used through.

\section{Heisenberg doubles}\label{chapter2}

In this section we shortly describe the basic notions concerning $\mathbb{Z}_2$-graded Heisenberg doubles and sketch how the Heisenberg double construction works. The exposition is structured in a way that is similar to \cite{Ka3}, it is however we generalised it such that it works in the $\mathbb{Z}_2$-graded setting.

We start from a short description of $\mathbb{Z}_2$-graded Hopf algebras, and using the Hopf action we will define a smash product of a Hopf algebra with its dual. We also discuss the connection between the Heisenberg double and the handle algebras of Chern-Simons theory as constructed using combinatorial quantisation due to Alekseev-Schomerus \cite{Alekseev:1994pa,Alekseev:1994au,Alekseev:1995rn, Alekseev:1996jz, Alekseev:1996ns,VS} for the non-graded gauge groups and \cite{AGPS1, AGPS2} for the graded, (possibly) non-semisimple gauge groups.

Let us consider a $\mathbb{Z}_2$-graded Hopf-algebra $(\mathcal{A},m,\eta,\Delta,\epsilon,\gamma)$, where $\mathcal{A}$ is a $\mathbb{Z}_2$-graded vector space equipped with the multiplication $m:\mathcal{A}\otimes\mathcal{A}\to\mathcal{A}$, an unit $\eta:\mathbb{C}\to \cA$, a co-multiplication $\Delta: \cA\to\cA\otimes\cA$, a co-unit $\epsilon:\cA\to\mathbb{C}$ and an antipode $\gamma:\cA\to\cA$. $\cA$ decomposes into a direct sum of two sub-spaces $\cA = \cA_0 \oplus \cA_1$, which are called even and odd respectively. We denote the degree of an element $x\in\cA_i$ by $|x|=i$, and we will call an element $x$ even if $|x|=0$ and odd otherwise. The graded tensor product of two algebras $\cA$ and $\mathcal{B}$ is then defined by the following equation for $a_1, a_2\in\cA$, $b_1,b_2\in \mathcal{B}$,
\begin{align}
	(a_1 \otimes b_1)\cdot(a_2 \otimes b_2) = (-1)^{|b_1||a_2|} a_1 a_2 \otimes b_1 b_2.
\end{align}
 The maps $m$ and $\eta$ are subjected to the associativity and unitality relations
\begin{align}\label{hopf_alg_axioms_1}
&	m \circ (m\otimes id) = m \circ (id \otimes m),\\
&	m \circ (\eta \otimes id) = id = m \circ (id \otimes \eta) ,
\end{align}
while maps $\Delta$ and $\epsilon$ have to satisfy the co-associativity and co-unitality relations
\begin{align}
&	(\Delta\otimes id)\circ \Delta = \Delta\circ (id \otimes \Delta),\\
&	(\epsilon \otimes id)\circ \Delta = id = (id \otimes \epsilon)\circ\Delta.
\end{align}
Moreover, the co-product $\Delta$ and co-unit $\epsilon$ are algebra homomorphisms, and the antipode $\gamma$ is a graded algebra anti-homomorphism and a graded co-algebra anti-homomorphism which satisfies the relations
\begin{align}\label{hopf_alg_axioms_3}
	m \circ (id\otimes \gamma)\circ \Delta = m\circ (\gamma\otimes id)\circ \Delta = \eta\circ\epsilon .
\end{align}
All the above maps are grade preserving. 

Moreover, we consider a graded Hopf algebra $(\mathcal{A}^*,{\hat m},{\hat \eta},{\hat \Delta},{\hat \epsilon},{\hat \gamma})$ which is dual to $\cA$. The Hopf algebras $\mathcal{A}$ and $\mathcal{A}^*$ are dual in a sense that the vector spaces $\mathcal{A}$ and $\mathcal{A}^*$ are dual as vector spaces, and there exists a non-degenerate duality pairing (also called a \textit{Hopf pairing}) $(,): \mathcal{A}\times\mathcal{A}^*\rightarrow \mathbb{C}$, for which the following relations are satisfied
\begin{gather}\label{duality_bracket_eq1}
 (x, f g) = (\Delta(x),f\otimes g), \qquad (x y, f) = (x\otimes y, {\hat \Delta}(f)), 
\end{gather}
between multiplications and co-multiplications,
\begin{gather}
 (\eta(1) , f) = {\hat \epsilon}(f), \qquad \epsilon(x) = (x,{\hat \eta}(1)), 
\end{gather}
between unit and co-unit maps,
\begin{gather}\label{duality_bracket_eq3}
 (\gamma(x),f) = (x,{\hat \gamma}(f)) ,
\end{gather}
and between antipodes, where
\begin{align*}
 (x\otimes y, f\otimes g) = (-1)^{|y||f|} (x,f)(y,g) ,
\end{align*}
for $x,y\in\mathcal{A}$, $f,g\in\mathcal{A}^*$.

The ordinary tensor product $\cA^* \otimes \cA$ has a straight-forward product given by $(1\otimes x)(f\otimes 1) = f\otimes x $. 
However, in order to construct a Heisenberg double, we are interested in equipping the space $\cA^* \otimes \cA$ with a non-trivial algebra structure between the elements which belong to the subalgebras $\cA$ and $\cA^*$. In order to achieve that, we will use the Hopf pairing $(,)$ to define a left action of $\cA$ on $\cA^*$ and consequently a smash product algebra $\mathcal{A}^* \rtimes \mathcal{A}$. \\\vspace{-5pt}

Using the duality pairing $(,)$ we can define a left action $\triangleright$ of a graded Hopf-algebra $\mathcal{A}$ on $\mathcal{A}^*$ given by
\begin{equation}\label{action_of_cA_on_cA*}
 x\otimes f \mapsto \sum_{(f)} (-1)^{|f_{(1)}||f_{(2)}|} (x, f_{(2)}) f_{(1)} =: x \triangleright f ,
\end{equation}
where $x\in\cA, f\in\cA^*$ and where we denote the co-product ${\hat \Delta}(f) = \sum_{(f)} f_{(1)}\otimes f_{(2)}$ using the usual Sweedler notation \cite{Sw}. The action~\eqref{action_of_cA_on_cA*} makes $\cA^*$ into a module algebra over the Hopf algebra $\cA$, i.e.\
the action is compatible with the multiplication in $\cA^*$ in the sense that
\begin{equation}\label{eq:x-fg}
 x \triangleright (f g) = \sum_{(x)} (-1)^{|f| |x_{(2)}|} (x_{(1)} \triangleright f) \, (x_{(2)} \triangleright g) ,
\end{equation}
where $x\in\cA$ and $f,g\in\cA^*$. Using a left action $\triangleright$ one can construct a smash product algebra $H(\cA) = \cA^* \rtimes \cA$ by defining the multiplication
\begin{equation}\label{heisenbergmultiplication}
(f\otimes x)(g\otimes y) = \sum_{(x)} (-1)^{|g||x_{(2)}|} f (x_{(1)} \triangleright g)\otimes x_{(2)} y ,
\end{equation}
where $x,y\in \mathcal{A}$, $f,g\in \mathcal{A}^*$. 

\begin{defn}\label{sec02:heisenberg-double}
A {\normalfont{Heisenberg double}} of a graded Hopf algebra $\cA$ is the smash product algebra $H(\cA) = \cA^* \rtimes \cA$ with the multiplication given by the equation \eqref{heisenbergmultiplication}.
\end{defn}

The Heisenberg double has $\mathcal{A}$ and $\mathcal{A}^*$ as subalgebras through canonical embeddings $\mathcal{A}^*\ni f\mapsto f\otimes1\in H(\mathcal{A})$ and $\mathcal{A}\ni x \mapsto 1\otimes x \in H(\mathcal{A})$. 

It is important to note that the Heisenberg double $H(\cA)$ is not a Hopf algebra. The algebra structure given by $\eqref{heisenbergmultiplication}$ is not compatible with the co-products on $\Delta$, ${\hat \Delta}$ defined on the initial Hopf algebras $\cA$ and $\cA^*$. In this it differs from the Drinfeld double algebra (which will be discussed in Section \ref{chapter3}), which is a (quasi-triangular) Hopf algebra and not only an algebra.

It will be convenient to recast the definition above in a basis dependent way. In order to do that, we first choose a basis of $\cA$. The basis will be given by a collection of vectors $\{e_\alpha\}_{\alpha\in I}$, where $I$ is a set. Then, the multiplication and co-multiplication of the basis elements is given by
\begin{align}
& e_\alpha e_\beta = \sum_{\gamma\in I} m^\gamma_{\alpha\beta} e_\gamma,
&& \Delta(e_\alpha) = \sum_{\beta, \gamma\in I} \mu^{\beta\gamma}_\alpha e_\beta\otimes e_\gamma ,
\end{align}
where $m^\gamma_{\alpha\beta}$ and $\mu^{\beta\gamma}_\alpha$ are respectively multiplication and co-multiplication coefficients. With the choice of a basis $\{e^\alpha\}_{\alpha\in I}$ of $\cA^*$ dual to $\{e_\alpha\}_{\alpha\in I}$ in the sense
\begin{equation}\label{duality}
(e_\alpha,e^\beta) = \delta_\alpha^\beta ,
\end{equation}
the multiplication and co-multiplication on $\cA^*$ are as follows
\begin{align}
&e^\alpha e^\beta = \sum_{\gamma\in I} (-1)^{|\alpha||\beta|} \mu_\gamma^{\alpha\beta} e^\gamma,
&&{\hat \Delta}(e^\alpha) = \sum_{\beta, \gamma\in I} (-1)^{|\beta||\gamma|} m_{\beta\gamma}^\alpha e^\beta\otimes e^\gamma .
\end{align}
The signs in the expressions above are the consequence of the properties of duality bracket. The Heisenberg double $H(\mathcal{A})$ is thus spanned by a collection of basis elements $\{ e^\alpha \otimes e_\beta\}_{\alpha,\beta\in I}$, written in terms of the basis elements of $\cA$ and $\cA^*$. The multiplication \eqref{heisenbergmultiplication} on basis elements $e^\alpha \otimes e_\beta$ has the following form
\begin{equation}
 (e^\alpha \otimes e_\beta)(e^\gamma \otimes e_\delta) = \sum_{\epsilon,\pi,\rho,\sigma,\tau\in I} (-1)^{|\beta||\gamma| + |\pi||\epsilon| + |\pi||\alpha| + |\epsilon|} m^\gamma_{\pi\epsilon} m^\tau_{\rho\delta} \mu^{\epsilon\rho}_\beta \mu^{\alpha\pi}_\sigma e^\sigma \otimes e_\tau.
\end{equation} 

Although, as we have mentioned previously, the Heisenberg double is not a Hopf algebra, one is interested in a way to encode the co-products on the initial graded Hopf algebras $\mathcal{A}$ and $\mathcal{A}^*$.

\begin{defn}
 Let us define an element called the {\normalfont{canonical element}} $W \in H(\mathcal{A})\otimes H(\mathcal{A})$ that has the following properties
\begin{align}\begin{aligned}
 W^{-1} \big( (1\otimes1)\otimes (1\otimes e_\alpha) \big) W &= \sum_{\beta, \gamma\in I} \mu^{\beta\gamma}_\alpha (1 \otimes e_\beta) \otimes (1 \otimes e_\gamma), \\
W\big((e^\alpha \otimes 1) \otimes (1\otimes1)\big)W^{-1} &= \sum_{\beta, \gamma\in I} (-1)^{|\beta||\gamma|} m_{\beta\gamma}^\alpha (e^\beta \otimes 1) \otimes (e^\gamma \otimes 1) .
\end{aligned}\end{align} 
\end{defn}

Given the above definition, one can show that the following two propositions about the canonical element $W$ are true.

\begin{prop}
 The canonical element $W$ defined above satisfies the {\normalfont{graded pentagon relation}}
\begin{equation}\label{pentagonrelation} 
W_{12}W_{13}W_{23} = W_{23}W_{12},
\end{equation}
where we use a notation for which $W_{12} = W \otimes (1\otimes1)$, $W_{23} = (1\otimes1)\otimes W$ and $W_{13} = \sigma_{12} [(1\otimes1) \otimes W]$, where $\sigma$ is an exchange operator.
\end{prop}

\begin{prop}
The canonical element $W$ can be expressed in terms of basis elements in the following way
\begin{equation}\label{canonicalelementdef}
W = \sum_{\alpha\in I} (-1)^{|\alpha|} (1\otimes e_\alpha)\otimes (e^\alpha \otimes 1) ,
\end{equation} 
\end{prop}

\begin{re}
 From Definition \ref{sec02:heisenberg-double} follows that the Heisenberg double $H(\cA^*)$ is spanned by a collection of elements $\{{\tilde e}_\alpha \otimes {\tilde e}^\beta \}_{\alpha,\beta\in I}$ subjected to the following relations
\begin{equation}\label{sec02:eq-crossing-relations-heisenberg-dual}({\tilde e}_\alpha \otimes \tilde e^\beta)(\tilde e_\gamma \otimes {\tilde e}^\delta) = \sum_{\epsilon,\pi,\rho,\sigma,\tau\in I} (-1)^{|\rho||\pi|+|\rho||\epsilon|+|\pi||\delta|} \mu^{\rho\epsilon}_\gamma \mu^{\pi\delta}_\tau m^\beta_{\epsilon\pi} m^\sigma_{\alpha\rho} \tilde e_\sigma \otimes \tilde e^\tau .
\end{equation}
We will denote the canonical element of this Heisenberg double as ${\tilde W} = \sum_{\alpha\in I} ({\tilde e}_\alpha \otimes 1) \otimes (1 \otimes {\tilde e}^\alpha )$. It satisfies a ``reversed'' pentagon equation of the form
\begin{align}
 {\tilde W}_{12} {\tilde W}_{23} = {\tilde W}_{23} {\tilde W}_{13} {\tilde W}_{12} .
\end{align}
\end{re}
We can relate those 2 algebras by the means of the following proposition:
\begin{prop}\label{proposition_heis_dual_heis_map}
 There exists an algebra anti-isomorphism $\xi: H(\cA^*) \to H(\cA)$ given by
 \begin{align}
  \xi({\tilde e}_\alpha \otimes 1) = \sum_{\beta\in I} \gamma_\alpha^\beta (1\otimes e_\beta), && \xi(1\otimes {\tilde e}^\alpha) = \sum_{\beta\in I} (-1)^{|\alpha|} (\gamma^{-1})^\alpha_\beta (e^\beta \otimes 1).
 \end{align}
 Since $ H(\cA^*) \ni x \mapsto (-1)^{|x|} x \in H(\cA^*)$ is a Heisenberg double isomorphism, the map $ H(\cA^*) \ni x \mapsto (-1)^{|x|} \xi(x) \in H(\cA)$ is also an algebra anti-isomorphism.
\end{prop}
The anti-isomorphism can be implemented on representation spaces in terms of super-transposition (i.e. the graded analogue of ordinary transposition). The super-transposition for square even $(n|m)$-matrices, i.e. linear transformations belonging to the space of $\Hom(\bC^{n|m},\bC^{n|m})$, is given by 
\begin{equation}
 \left( \begin{array}{cc} A & B\\C& D \end{array} \right)^\text{st} = \left( \begin{array}{cc} A^\text{t} & C^\text{t}\\-B^\text{t}& D^\text{t} \end{array} \right) ,
\end{equation}
where $A\in \Hom(\mathbb{C}^n,\mathbb{C}^n)$, $D\in \Hom(\mathbb{C}^m,\mathbb{C}^m)$ are even, $B\in \Hom(\mathbb{C}^m,\mathbb{C}^n)$, $C\in \Hom(\mathbb{C}^n,\mathbb{C}^m)$ are odd, and $^\text{t}$ denotes an ordinary, not graded matrix transposition.

\begin{re}
 To keep the notation compact, from now on we will denote the elements $1\otimes e_\alpha $ and $e^\alpha\otimes1 $ of the Heisenberg double $H(\cA)$ simply as $e_\alpha$ and $e^\alpha$ respectively.
\end{re}

Finally, we would like to make a connection between the Heisenberg double and the handle algebra of the Chern-Simons theory. While the Heisenberg double is defined for arbitrary Hopf algebras, for the purposes of investigating this connection we will restrict ourselves to $\cA$ being a quasi-triangularity Hopf algebra, i.e. that it does admit an $R$-matrix. As mentioned in the introduction, the handle algebra arises as a graph algebra associated to the torus $\Sigma_{1,0}$ in the context of the Chern-Simons theory with a gauge group $G$. To be precise, it can be defined in the following way.

\begin{defn}
 Let $\cA$ be a quasi-triangular $\mathbb{Z}_2$-graded Hopf algebra. The {\normalfont{handle algebra}} $\cT(\cA)$ associated to the graded Hopf algebra $\cA$ is defined as an algebra spanned by the elements $f_a^i f_b^i$ such that the universal elements $A,B\in\cA\otimes \cT(\cA)$
 \begin{equation}
  A = \sum_{i\in I} e_i \otimes f_a^i , \qquad B = \sum_{i\in I} e_i \otimes f_b^i ,
 \end{equation}
 satisfy the following quadratic equations
 \begin{gather}\label{ch02:handle-algebra-eq-1}
  A_{13} R_{12} A_{23} = R_{12} (\Delta\otimes id)A, \\
  B_{13} R_{12} B_{23} = R_{12} (\Delta\otimes id)B, \label{ch02:handle-algebra-eq-2} \\
  R_{12}^{-1} A_{13} R_{12} B_{23} = B_{23} R'_{12} A_{13} R_{12} .\label{ch02:handle-algebra-eq-3}
 \end{gather} 
\end{defn}
While in the definition of the handle algebra we start from universal elements, the Heisenberg double $H(\cA^*)$ also admits a description using its own universal elements, which are defined as follows.

\begin{defn}\label{ch02:def-universal-element-G}
 The universal element $G \in \cA\otimes \cA^*\subset \cA \otimes H(\cA^*)$ is defined by the equation
 \begin{equation}\label{universalG}
   G = \sum_{i\in I} e_i \otimes e^i,
 \end{equation}
 where $e_i\in\cA$ and $e^i\in H(\cA^*)$.
\end{defn}

\begin{defn}\label{ch02:def-universal-element-N}
 The universal elements $N_\pm \in \cA\otimes \cA \subset \cA \otimes H(\cA^*)$ defined by
 \begin{equation}
  N_+ = (id\otimes \iota) R' \qquad N_- = (id\otimes \iota) R^{-1},
 \end{equation}
 where $\iota: \cA\to H(\cA^*)$ is the canonical embedding of $\cA$ into $H(\cA^*)$.
\end{defn}

Then, one can define an isomorphism between those two algebras on the level of universal elements, which descends directly to an isomorphism on the elements.

\begin{prop}\label{ch02:prop-heisenberg-handle-isomorphism}
 There exists an isomorphism $\varphi: \cT(\cA) \to H(\cA^*)$ defined by the equations
 \begin{equation}
   (id\otimes \varphi) A = N_+ G N_-^{-1} , \qquad (id\otimes \varphi) B = N_+ N_-^{-1} ,
 \end{equation}
 between the handle algebra of the $\mathbb{Z}_2$-graded Hopf algebra $\cA$ and the Heisenberg double $H(\cA^*)$.
\end{prop}

As the proof of this assertion was missing in the literature regarding even the non-graded case, we will give it below, after producing some additional propositions and lemmas.

\begin{prop}\label{appendixA:prop1}
 The universal element $G \in \cA\otimes \cA^*\subset \cA \otimes H(\cA^*)$ defined in Definition \ref{ch02:def-universal-element-G} satisfies the equation
 \begin{equation}\label{appendixA:eq-universal-element-G-relation}
  (\Delta \otimes id) G = G_{13} G_{23} .
 \end{equation}
\begin{proof}
 Using the co-product for the elements $e_i$ and the product for the elements $e^i$ one gets
 \begin{align*}
  (\Delta \otimes id) G &= \sum_{i\in I} \Delta(e_i) \otimes e^i = \sum_{i,j,k\in I} e_j \otimes e_k \otimes \mu_i^{jk} e^i = \sum_{j,k\in I} (-1)^{|j||k|} e_j \otimes e_k \otimes e^j e^k = \\
  &= \sum_{j\in I} (e_j \otimes 1 \otimes e^j) \sum_{k\in I} (1 \otimes e_k \otimes e^k) = G_{13} G_{23} .
 \end{align*}
\end{proof}
\end{prop}

\begin{prop}\label{appendixA:prop2}
 The relations \eqref{sec02:eq-crossing-relations-heisenberg-dual} are equivalent to the requirement that for all $\xi\in\cA$ the universal element $G \in \cA\otimes \cA^* \subset \cA \otimes H(\cA^*)$ ought to satisfy
 \begin{equation}
  G (1\otimes \iota(\xi)) = (id\otimes \iota)\Delta^{op}(\xi) G,
 \end{equation}
 where $\iota: \cA \to H(\cA^*)$ is a canonical embedding of $\cA$ in $H(\cA^*)$.
\begin{proof}
 Setting $\xi = e_\gamma$ and dropping $\iota$ for clarity one can calculate using the decomposition of the element $G$
 \begin{align*}
  \text{LHS} &= \sum_{\beta\in I} e_\beta \otimes e^\beta e_\gamma , \\
  \text{RHS} &= \sum_{\epsilon,\sigma\in I} (-1)^{|\sigma||\epsilon|} \mu_\gamma^{\sigma\epsilon} (e_\epsilon \otimes e_\sigma) (e_\rho \otimes e^\rho) = \sum_{\epsilon,\sigma\in I} (-1)^{|\sigma|(|\epsilon|+|\rho|)} \mu_\gamma^{\sigma\epsilon} e_\epsilon e_\rho \otimes e_\sigma e^\rho = \\
  &= \sum_{\beta,\epsilon,\sigma\in I} (-1)^{|\sigma|(|\epsilon|+|\rho|)}  \mu_\gamma^{\sigma\epsilon} m_{\epsilon\rho}^\beta e_\beta \otimes e_\sigma e^\rho =  \sum_{\beta\in I} e_\beta \otimes \left[ \sum_{\epsilon,\sigma\in I} (-1)^{|\sigma|(|\epsilon|+|\rho|)}  \mu_\gamma^{\sigma\epsilon} m_{\epsilon\rho}^\beta e_\sigma e^\rho \right] .
 \end{align*}
 Because the basis elements $e_\beta$ are linearly independent, the elements in the second tensor factor have to be equal for each $\beta$ independently, which leads to the relation
 \begin{equation*}
  e^\beta e_\gamma = \sum_{\epsilon,\sigma\in I} (-1)^{|\sigma|(|\epsilon|+|\rho|)}  \mu_\gamma^{\sigma\epsilon} m_{\epsilon\rho}^\beta e_\sigma e^\rho ,
 \end{equation*}
 which was to be shown.
\end{proof}
\end{prop}

\begin{lem}\label{appendixA:lemma-quasitriangular-N}
 The universal elements $N_+$, $N_-$ given in Definition \ref{ch02:def-universal-element-N} satisfy
 \begin{equation}
  (\Delta\otimes id) (N_\pm) = (N_\pm)_{23} (N_\pm)_{13} .
 \end{equation}
 \begin{proof}
  This follows directly from the quasi-triangularity relations \eqref{ch3:eq-quasi-triangularity} for the $R$-matrix.
 \end{proof}
\end{lem}

\begin{lem}\label{appendixA:lemma-quasitriangular-N-2}
 The universal elements $N_+$, $N_-$ satisfy
 \begin{equation}
  R_{12} (N_\pm)_{23} (N_\pm)_{13} = (N_\pm)_{13} (N_\pm)_{23} R_{12}, \qquad R_{12} (N_+)_{23} (N_-)_{13} = (N_-)_{13} (N_+)_{23} R_{12} .
 \end{equation}
 \begin{proof}
  The first equation follows directly from Lemma \ref{appendixA:lemma-quasitriangular-N} and the fact that the $R$-matrix is an intertwiner \eqref{ch3:eq-intertwining-R}. The second one can be obtained from the quasi-triangularity \eqref{ch3:eq-quasi-triangularity} and the intertwining property of the $R$-matrix.
 \end{proof}
\end{lem}

\begin{lem}\label{appendixA:quadratic-eq-for-N-lemma}
 The universal element $N = N_+ N_-^{-1}$ satisfies the following relation
 \begin{equation}
  R(\Delta\otimes id) N = N_{13} R_{12} N_{23} .
 \end{equation}
 \begin{proof}
  Using above lemmas, one can show that
  \begin{align*}
   (\Delta\otimes id)N &= (N_+)_{23} (N_+)_{13} (N_-)^{-1}_{13} (N_+)^{-1}_{23} = R^{-1}_{12} (N_+)_{13} (N_+)_{23} R_{12} (N_-)^{-1}_{13} (N_+)^{-1}_{23} = \\
   &= R^{-1}_{12} (N_+)_{13} (N_-)^{-1}_{13} R_{12} (N_+)_{23} (N_+)^{-1}_{23} = R^{-1}_{12} N_{13} R_{12} N_{23} ,
  \end{align*}
  which completes the proof.
 \end{proof}
\end{lem}

\begin{prop}\label{propositionA6}
 The universal elements $N_+$, $N_-$ and $G$ satisfy the following quadratic equations
 \begin{equation}
  (N_+)_{13} R'_{12} G_{23} = G_{23} (N_+)_{13}, \quad (N_-)_{13} R^{-1}_{12} G_{23} = G_{23} (N_-)_{13},
 \end{equation}
 \begin{proof}
  Let us assume that the $R$ matrix decomposes into elements $f_i, g_i \in \cA$ as follows
  \begin{equation*}
   R = \sum_{i\in I} f_i \otimes g_i, \qquad R' = \sum_{i\in I} (-1)^{|f_i|} g_i \otimes f_i .
  \end{equation*}
  Then, using the relation
  \begin{equation*}
   (id \otimes \Delta^{op})(R') = R'_{13} R'_{12} ,
  \end{equation*}
  which follows from the left equation \eqref{ch3:eq-quasi-triangularity}, as well as using Proposition \ref{appendixA:prop2} one can prove that
  \begin{align*}
   G_{23} (N_+)_{13} &= \sum_{i\in I} (-1)^{|f_i|} (1 \otimes G) (g_i \otimes \iota(f_i)) = \sum_{i\in I} (-1)^{|f_i|} (g_i \otimes G) (1\otimes1 \otimes \iota(f_i)) = \\
   &= \sum_{i\in I} (-1)^{|f_i|} (g_i \otimes (id\otimes\iota) \Delta^{op}(f_i)) G_{23} = \\
   &= \Big[(id\otimes id\otimes\iota) (id \otimes \Delta^{op})(R') \Big] G_{23} = \Big[(id\otimes id\otimes\iota) R'_{13} R'_{12} \Big] G_{23} = \\
   &= (N_+)_{13} R'_{12} G_{23} .
  \end{align*}
 The second equation, with the use of the identity
  \begin{equation*}
   (id \otimes \Delta^{op})(R^{-1}) = R^{-1}_{13} R^{-1}_{12} ,
  \end{equation*}
 which follows from the right equation \eqref{ch3:eq-quasi-triangularity} can be show analogously.
 \end{proof}
\end{prop}

Finally, we can prove Proposition \ref{ch02:prop-heisenberg-handle-isomorphism}.

\begin{proof}[Proof (of Proposition \ref{ch02:prop-heisenberg-handle-isomorphism})]
 In order to do that, we need to show that the isomorphism $\varphi$ preserves the defining relations \eqref{ch02:handle-algebra-eq-1}-\eqref{ch02:handle-algebra-eq-3}. For the equation \eqref{ch02:handle-algebra-eq-1} we can perform the following chain of application of identities given above
 \begin{align*}
  &R_{12} (\Delta \otimes \varphi) A = R_{12} (\Delta \otimes id)(N_+ G N_-^{-1}) = R_{12} \Big[ \underbrace{(\Delta \otimes id)N_+}_{=(N_+)_{23} (N_+)_{13}}\Big] \Big[ \underbrace{(\Delta \otimes id)G}_{=G_{13}G_{23}}\Big] \Big[\underbrace{(\Delta \otimes id)N_-^{-1}}_{=(N_-)^{-1}_{13} (N_-)^{-1}_{23}}\Big] = \\
  &= \underbrace{R_{12} (N_+)_{23} (N_+)_{13}}_{=(N_+)_{13} (N_+)_{23} R_{12}} G_{13}G_{23} (N_-)^{-1}_{13} (N_-)^{-1}_{23} = (N_+)_{13} (N_+)_{23} R_{12} G_{13} \underbrace{G_{23} (N_-)^{-1}_{13}}_{=R_{12}(N_-)^{-1}_{13}G_{23}} (N_-)^{-1}_{23} = \\
  &= (N_+)_{13} \underbrace{(N_+)_{23} R_{12} G_{13}}_{=G_{13}(N_+)_{23}} R_{12}(N_-)^{-1}_{13}G_{23} (N_-)^{-1}_{23} = (N_+)_{13} G_{13} \underbrace{(N_+)_{23} R_{12}(N_-)^{-1}_{13}}_{=(N_-)^{-1}_{13} R_{12} (N_+)_{23}} G_{23} (N_-)^{-1}_{23} = \\
  &= (N_+)_{13} G_{13} (N_-)^{-1}_{13} R_{12} (N_+)_{23} G_{23} (N_-)^{-1}_{23} = (id\otimes id\otimes \varphi) A_{13} R_{12} A_{23} .
 \end{align*}
 The equation \eqref{ch02:handle-algebra-eq-2} is a direct consequence of Lemma \ref{appendixA:quadratic-eq-for-N-lemma}. Finally, the equation \eqref{ch02:handle-algebra-eq-3} is shown using the following chain of identities
 \begin{align*}
  &(id\otimes id\otimes \varphi) R^{-1}_{12} A_{13} R_{12} B_{23} = R^{-1}_{12} (N_+)_{13} G_{13} \underbrace{(N_-)^{-1}_{13} R_{12} (N_+)_{23}}_{=(N_+)_{23}R_{12}(N_-)^{-1}_{13}} (N_-)^{-1}_{23} = \\
  &= R^{-1}_{12} (N_+)_{13} \underbrace{G_{13} (N_+)_{23}}_{=(N_+)_{23}R_{12}G_{13}} \underbrace{R_{12}(N_-)^{-1}_{13} (N_-)^{-1}_{23}}_{=(N_-)^{-1}_{23}(N_-)^{-1}_{13}R_{12}} =\\
  &= \underbrace{R^{-1}_{12} (N_+)_{13} (N_+)_{23}}_{=(N_+)_{23} (N_+)_{13}R^{-1}_{12}} R_{12} \underbrace{G_{13} (N_-)^{-1}_{23}}_{=R'_{12} (N_-)^{-1}_{23} G_{13}} (N_-)^{-1}_{13}R_{12} =\\
  &= (N_+)_{23} \underbrace{(N_+)_{13} R'_{12} (N_-)^{-1}_{23}}_{=(N_-)^{-1}_{23}R'_{12}(N_+)_{13}} G_{13} (N_-)^{-1}_{13}R_{12} = (N_+)_{23} (N_-)^{-1}_{23} R'_{12} (N_+)_{13} G_{13} (N_-)^{-1}_{13} R_{12} = \\
  &= (id\otimes id\otimes \varphi) B_{23} R'_{12} A_{13} R_{12} .
 \end{align*}
 This completes the proof of the proposition.
\end{proof}

\section{Drinfeld doubles}\label{chapter3}

In this section, we will present the definition of the Drinfeld double $D(\cA)$ of a $\mathbb{Z}_2$-graded Hopf algebra $\cA$ and remind ourselves some facts about the universal element $R$ satisfying the Yang-Baxter equation. Then, we will describe an algebra morphism between $H(\cA)\otimes H(\cA^*)$ and $D(\cA)$, which constitutes a $\mathbb{Z}_2$-graded generalisation of a morphism described in \cite{Ka3}. Furthermore, we state the relation between the universal elements of the Heisenberg doubles and the universal $R$-matrix of the Drinfeld double. We also discuss the connection between the Drinfeld double and the loop algebra $\mathcal{L}$ of Chern-Simons theory and its semi-direct product $\bar{\mathcal{L}}$ as constructed using combinatorial quantisation due to Alekseev-Schomerus \cite{Alekseev:1994pa,Alekseev:1994au,Alekseev:1995rn,Alekseev:1996ns,VS} for the non-graded gauge groups and \cite{AGPS1, AGPS2} for the $\mathbb{Z}_2$-graded, (possibly) non-semisimple gauge groups. \\

Let us consider again a $\mathbb{Z}_2$-graded Hopf algebra $(\mathcal{A},m,\eta,\Delta,\epsilon,\gamma)$, subjected to the axioms \eqref{hopf_alg_axioms_1}-\eqref{hopf_alg_axioms_3}. With the choice of a basis $\{e_\alpha\}_{\alpha\in I}$ which algebraically spans $\cA$ we can describe the multiplication, co-multiplication and an antipode 
\begin{align}
& e_\alpha e_\beta = \sum_{\gamma\in I} m^\gamma_{\alpha\beta} e_\gamma,
&& \Delta(e_\alpha) = \sum_{\beta, \gamma\in I} \mu^{\beta\gamma}_\alpha e_\beta\otimes e_\gamma ,
&& \gamma(e_\alpha) = \sum_{\beta\in I} \gamma_\alpha^\beta e_\beta .
\end{align}
In addition to the algebra $\cA$ we can consider the algebra $\cA^*$, that is a Hopf algebra dual to $\cA$. In terms of a basis $\{e^\alpha\}_{\alpha\in I}$ of $\cA^*$ the multiplication, co-multiplication and an antipode relations of this Hopf algebra are as follows
\begin{align}
&e^\alpha e^\beta = \sum_{\gamma\in I} (-1)^{|\alpha||\beta|} \mu_\gamma^{\alpha\beta} e^\gamma,
&&{\hat \Delta}(e^\alpha) = \sum_{\beta, \gamma\in I} (-1)^{|\beta||\gamma|} m_{\beta \gamma}^\alpha e^\beta\otimes e^\gamma,
&&{\hat \gamma}(e^\alpha) = \sum_{\beta\in I} \gamma^\alpha_\beta e^\beta.
\end{align}
They are dual to each other with respect to a duality bracket $(,)$ satisfying the relations \eqref{duality_bracket_eq1}-\eqref{duality_bracket_eq3}, and given explicitly on the basis by \eqref{duality}.

Given those two graded Hopf algebras, it is possible to define a quasi-triangular Hopf algebra as a double cross product of Hopf algebras (e.g. \cite{Majidbook} for non-graded case). The definition for the $\mathbb{Z}_2$-graded case is a straightforward generalisation of the non-graded Drinfeld double. Explicitly, we can consider a Hopf algebra $D(\cA)=\{ x\otimes f | x\in\cA, f\in\cA^*\}$ equipped with the product 
\begin{equation}\label{drinfeld_definition_eq1}\begin{aligned}
 (x\otimes f)(y\otimes g) &= \sum_{(y), (f)} (-1)^{|y_{(1)}| |f| + |y_{(2)}| (|f_{(1)}|+|f_{(2)}|) + |y_{(3)}| |f_{(1)}| + |f_{(2)}| |f_{(3)}| } \times \\
 &\qquad\times (-1)^{|f_{(1)}| (|f_{(2)}|+|f_{(3)}|) + |x| (|y_{(1)}| + |f_{(3)}|) + |g| (|f_{(1)}| + |x_{(3)}|)} \times \\
 &\qquad\times ( y_{(1)}, \hat{\gamma}^{-1}(f_{(3)}) ) ( y_{(3)}, f_{(1)} ) x y_{(2)} \otimes f_{(2)} g ,
\end{aligned}\end{equation}
the co-product
\begin{align}
 \Delta(x\otimes f) = \sum_{(x), (f)} (-1)^{|f_{(1)}||f_{(2)}|+|x_{(2)}||f_{(2)}|} x_{(1)} \otimes f_{(2)} \otimes x_{(2)} \otimes f_{(1)} , 
\end{align}
and the antipode
\begin{align}\label{drinfeld_definition_eq3}
 \gamma(x\otimes f) = (-1)^{|x||f|} (1 \otimes {\hat \gamma}^\text{cop}(f))(\gamma(x)\otimes 1),
\end{align}
where ${\hat \gamma}^\text{cop} = {\hat \gamma}^{-1}$ is the antipode of $(\mathcal{A}^*)^\text{cop}$, i.e. the antipode of the Hopf algebra dual to $\cA$ equipped with the opposite co-product. We also use the following notation for the co-product
\begin{align*}
 (\Delta\otimes id)\Delta(x) = \sum_{(x)} x_{(1)} \otimes x_{(2)} \otimes x_{(3)}, && ({\hat \Delta}\otimes id){\hat \Delta}(f) = \sum_{(f)} f_{(1)} \otimes f_{(2)} \otimes f_{(3)} .
\end{align*}

\begin{defn}
A {\normalfont{Drinfeld double}} of a graded Hopf algebra $\cA$ is a quasi-triangular Hopf algebra $D(\cA)$ with the multiplication, co-multiplication and antipode given by the equations \eqref{drinfeld_definition_eq1}-\eqref{drinfeld_definition_eq3}.
\end{defn}

The Hopf algebra $D(\cA)$ has $\mathcal{A}$ and $(\mathcal{A}^*)^\text{cop}$ as subalgebras through canonical embeddings $(\mathcal{A}^*)^\text{cop} \ni f\mapsto 1\otimes f\in D(\mathcal{A})$ and $\mathcal{A}\ni x \mapsto x\otimes 1 \in D(\mathcal{A})$. 

In terms of a basis, the Drinfeld double $D(\cA)$ is spanned by a collection of elements $\{e_\alpha \otimes e^\beta \}_{\alpha,\beta\in I}$ which satisfy the following multiplication relations
\begin{equation}\label{susyEE}
\begin{split}
& (e_\alpha \otimes1) (1\otimes e^\beta) = e_\alpha \otimes e^\beta ,\\
& (e_\alpha \otimes1) (e_\beta \otimes 1) = \sum_{\gamma\in I} m^\gamma_{\alpha\beta} (e_\gamma\otimes 1),\\
& (1\otimes e^\alpha) (1\otimes e^\beta) = \sum_{\gamma\in I} (-1)^{|\alpha||\beta|}\mu_\gamma^{\alpha\beta} (1\otimes e^\gamma), \\
& (1\otimes e^{\alpha} ) (e_{\beta} \otimes 1)= \sum_{\substack{\gamma,\delta,\epsilon,\\\mu,\nu,\rho,\sigma\in I}} (-1)^{|\mu|(|\sigma|+|\delta|)+|\alpha|+|\gamma|}  m_{\nu\gamma}^{\mu}m_{\mu\delta}^{\alpha} \mu_{\beta}^{\epsilon \rho}\mu_{\rho}^{\sigma \nu} ({\gamma}^{-1})^{\delta}_{\epsilon} e_{\sigma} \otimes e^{\gamma} ,
\end{split}
\end{equation}
and co-multiplication relations
\begin{align}
  \Delta(e_\alpha\otimes1) = \sum_{\beta, \gamma\in I} \mu^{\beta\gamma}_\alpha (e_\beta\otimes1)\otimes (e_\gamma\otimes1), &&   \Delta(1\otimes e^\alpha) = \sum_{\beta, \gamma\in I} m_{\gamma\beta}^\alpha (1\otimes e^\beta)\otimes (1\otimes e^\gamma) ,
\end{align}
and is equipped with the antipode
\begin{align}
  \gamma(e_\alpha\otimes 1) = \sum_{\beta\in I} \gamma_\alpha^\beta e_\beta\otimes 1, && \gamma(1\otimes e^\alpha) = \sum_{\beta\in I} {({\hat\gamma}^\text{cop})}_\alpha^\beta 1\otimes e^\beta ,
\end{align}
where $(\hat\gamma^\text{cop})^\alpha_\beta=(\gamma^{-1})^\alpha_\beta$, which follows directly from \eqref{drinfeld_definition_eq1}-\eqref{drinfeld_definition_eq3}. Equivalently, instead of the fourth exchange relation in \eqref{susyEE} one can use the crossing relation
\begin{align}\label{susytt}
&\sum_{\gamma,\rho,\sigma\in I} (-1)^{|\beta||\sigma|+|\rho|} \mu^{\sigma\gamma}_{\alpha} m^\beta_{\gamma\rho} e_\sigma \otimes e^\rho = \sum_{\gamma,\rho,\sigma\in I} (-1)^{|\rho||\gamma| + |\rho|} m^\beta_{\rho\gamma} \mu^{\gamma\sigma}_\alpha (1\otimes e^\rho) ( e_\sigma \otimes 1),
\end{align}
which indeed defines the same algebra.

\begin{defn}
In the case of Drinfeld double one has a {\normalfont{canonical element}} $R \in D(\mathcal{A})\otimes D(\mathcal{A})$ called the {\normalfont{universal $R$-matrix}}
\begin{equation}\label{rmatrix_def}
 R = \sum_{\alpha\in I} (-1)^{|\alpha|} (e_\alpha \otimes 1) \otimes (1 \otimes e^\alpha).
\end{equation} 
\end{defn}

\begin{prop}\label{ch03:prop-quasi-triangularity}
The universal $R$-matrix satisfies the {\normalfont{quasi-triangularity relations}}

  \begin{equation}\label{ch3:eq-quasi-triangularity}
   (\Delta \otimes id)(R) = R_{13} R_{23} , \qquad (id \otimes \Delta)(R) = R_{13} R_{12} ,
  \end{equation}
and it is an intertwiner between $\cA$ and $\cA^{cop}$, i.e.
  \begin{equation}\label{ch3:eq-intertwining-R}
   R \Delta(\xi) = \Delta^{op}(\xi) R,
  \end{equation}
 for $\xi\in D(\cA)$.
\end{prop}

Proposition \ref{ch03:prop-quasi-triangularity} implies the following proposition.

\begin{prop}
The universal $R$-matrix satisfies the {\normalfont{Yang-Baxter equation}}
\begin{equation}\label{yang_baxter_equation} 
 R_{12} R_{13} R_{23} = R_{23} R_{13} R_{12} .
\end{equation}
where we use a notation for which $R_{12} = R \otimes (1\otimes 1)$, $R_{23} = (1\otimes 1)\otimes R$ and $R_{13} = \sum_{\alpha\in I} (-1)^{|\alpha|} (e_\alpha \otimes 1) \otimes (1\otimes 1) \otimes ( 1\otimes e^\alpha)$. 
\end{prop}
Now we want to show that there exists an algebra embedding of the Drinfeld double $D(\cA)$ in a tensor square of Heisenberg doubles $H(\cA)$. The following proposition is true: 
\begin{prop} A map $\eta: D(\cA) \to H(\cA)\otimes H(\cA^*)$ defined as follows \label{drinfeld_heisenberg_homomorphism}
\begin{equation}\begin{aligned}
 \eta^{(a,b)}(e_\alpha\otimes1) &= \sum_{\beta,\gamma\in I} (-1)^{a |\beta| + b|\gamma|} \mu_\alpha^{\beta\gamma} (1\otimes e_\beta)\otimes(\tilde e_\gamma \otimes 1), \\
 \eta^{(a,b)}(1\otimes e^\alpha) &= \sum_{\beta,\gamma\in I}(-1)^{a |\beta| + (b+1)|\gamma|} m^\alpha_{\gamma\beta} (e^\beta\otimes 1)\otimes(1 \otimes \tilde e^\gamma), 
\end{aligned}\end{equation}
for the parameters $a,b=0,1$ (which choice corresponds to composition with $x \mapsto (-1)^{|x|} x$ for $x\in H(\cA)$ or $x\in H(\cA^*)$) is an algebra homomorphism. 
\end{prop}

\begin{cor}\label{sec3:corollary:drinfeld-double-algebra-homomorphism}
 From Propositions \ref{proposition_heis_dual_heis_map} and \ref{drinfeld_heisenberg_homomorphism} directly follows that one has an algebra homomorphism $\phi: D(\cA) \to H(\cA)\otimes H(\cA^{op,cop})$ defined as follows
 \begin{equation}\begin{aligned}
 \phi^{(a,b)}(e_\alpha\otimes1) &= \sum_{\beta,\gamma,\delta\in I} (-1)^{a |\beta| + b|\delta|} \mu_\alpha^{\beta\gamma} \gamma_\gamma^\delta (1\otimes e_\beta)\otimes (1\otimes e_\delta), \\
 \phi^{(a,b)}(1\otimes e^\alpha) &= \sum_{\beta,\gamma,\delta\in I}(-1)^{a |\beta| + b|\delta|} (\gamma^{-1})^\gamma_\delta m^\alpha_{\gamma\beta} (e^\beta\otimes 1)\otimes (e^\delta \otimes 1).
 \end{aligned}\end{equation}
\end{cor}

Using the morphism $\eta$ from Proposition \ref{drinfeld_heisenberg_homomorphism} in addition to the canonical element $W$ for the Heisenberg double $H(\cA)$ and the (flipped) canonical element ${\tilde W}$ for $H(\cA^*)$, one can define the following two elements $W' = \sum_{\alpha\in I} (-1)^{(a+b+1)|\alpha|} ({\tilde e}_\alpha \otimes 1) (1 \otimes e^\alpha)$ and $W'' = \sum_{\alpha\in I} (-1)^{(a+b) |\alpha|} (1 \otimes  e_\alpha) \otimes ({\tilde e}^\alpha \otimes 1)$.
It can be shown that they satisfy a set of 6 pentagon-like equations
\begin{align}
 W'_{12} W'_{13} W_{23} = W_{23} W'_{12}, && {\tilde W}_{12} W'_{23} = W'_{23} W'_{13} {\tilde W}_{12}, \nonumber \\
 W_{12} W''_{13} W''_{23} = W''_{23} W_{12}, && W''_{12} {\tilde W}_{23} = {\tilde W}_{23} W''_{13} W''_{12}, \\
 W'_{12} {\tilde W}_{13} W''_{23} = W''_{23} W'_{12}, && W''_{12} W'_{23} = W'_{23} W_{13} W''_{12}. \nonumber
\end{align}
Then, we claim that one can construct the R-matrix of the Drinfeld double $D(\cA)$ as follows
\begin{prop}\label{Rmatrix-fourS}
 Under an algebra map $\eta$ one has the following relation
 \begin{equation}
  (\eta^{(a,b)}\otimes \eta^{(a,b)}) R = W''_{14} W_{13} {\tilde W}_{24} W'_{23} .
 \end{equation}
\end{prop}

\begin{re}
 To keep the notation compact, from now on we will denote the elements $1\otimes e^\alpha $ and $e_\alpha\otimes1 $ of the Drinfeld double $D(\cA)$ simply as $e^\alpha$ and $e_\alpha$ respectively.
\end{re}

Lastly, we want to comment about the relationship between the Drinfeld double and the graph algebra $\bar{\mathcal{L}}$ present in the combinatorial Chern-Simons theory (as defined, e.g., in the Definition 5 in \cite{Alekseev:1996jz}). The isomorphism between those objects has been previously proven in the non-graded, semisimple setting \cite{Nill:1996dv, Alekseev:1996jz}. We generalise it to the graded setting below. \\

In the combinatorial Chern-Simons theory one can define a \textit{gauged loop algebra}, which is a semi-direct product of a \textit{loop algebra} $\mathcal{L}$ with the quasi-triangular Hopf algebra $\cA$, associated with 1-punctured sphere $\Sigma_{0,1}$, which we will denote by the symbol $\bar{\mathcal{L}}$.
\begin{defn}
 Let $\cA$ be a quasi-triangular $\mathbb{Z}_2$-graded Hopf algebra. The gauged loop algebra $\bar{\mathcal{L}}$ associated to the Hopf algebra $\cA$ is defined as an algebra spanned by the elements $f^i g^j$ such that the universal elements $M_0,M_1\in\cA\otimes \bar{\mathcal{L}}$
 \begin{equation}
  M_0 = \sum_{i\in I} e_i \otimes f^i , \qquad M_1 = \sum_{i\in I} e_i \otimes g^i ,
 \end{equation}
 satisfy the following quadratic equations
 \begin{gather}\label{ch03:loop-algebra-eq-1}
  (M_0)_{13} R_{12} (M_0)_{23} = R_{12} (\Delta\otimes id)M_0, \\
  (M_1){13} R_{12} (M_1)_{23} = R_{12} (\Delta\otimes id)M_1, \label{ch03:loop-algebra-eq-2} \\
  R_{12}^{-1} (M_1)_{13} R_{12} (M_0)_{23} = (M_0)_{23} R'_{12} (M_1)_{13} (R')^{-1}_{12} .\label{ch03:loop-algebra-eq-3}
 \end{gather} 
\end{defn}
\begin{re}
 We would like to comment briefly about the difference between the above definition and the Definition 5 from \cite{Alekseev:1996jz}, as well as the difference between the equation \eqref{ch03:loop-algebra-eq-3} and the equation \eqref{ch02:handle-algebra-eq-3}, which distinguishes the algebra $\bar{\mathcal{L}}$ from the handle algebra $\mathcal{T}(\cA)$. The components $f^i$ span a subalgebra in $\bar{\mathcal{L}}$ which is isomorphic to $\cA$ and therefore the algebra $\bar{\mathcal{L}}$ admits an embedding $\iota : \cA \to \text{span}\{ f^i \}_{i\in I} \subset \bar{\mathcal{L}}$. The equation \eqref{ch03:loop-algebra-eq-3} is a consequence of the fact that this algebra acts as an adjoint on its compliment, i.e. the following holds
 \begin{equation*}
  [(id\otimes\iota) \Delta(x)] M_1 = M_1 [(id\otimes\iota) \Delta(x)] ,
 \end{equation*}
 where $x \in \cA$. This explains how our definition is related to the one from \cite{Alekseev:1996jz}. Also, above is not the case for the handle algebra $\cT(\cA)$.
\end{re}

As in the case of Heisenberg double previously (i.e. Definitions \ref{ch02:def-universal-element-G} and \ref{ch02:def-universal-element-N}), for Drinfeld double one can define the universal elements $G$, $N_+$ and $N_-$.
\begin{defn}\label{ch03:def-universal-element-G}
 The universal element $G \in \cA\otimes \cA^*\subset \cA \otimes D(\cA)$ is defined by the equation
 \begin{equation}\label{ch03:universalG}
   G = \sum_{i\in I} (-1)^{|i|} e_i \otimes e^i,
 \end{equation}
 where $e_i\in\cA$ and $e^i\in D(\cA)$.
\end{defn}

\begin{defn}\label{ch03:def-universal-element-N}
 The universal elements $N_\pm \in \cA\otimes \cA \subset \cA \otimes D(\cA)$ defined by
 \begin{equation}
  N_+ = (id\otimes \iota) R' \qquad N_- = (id\otimes \iota) R^{-1},
 \end{equation}
 where $\iota: \cA\to D(\cA)$ is the canonical embedding of $\cA$ into $D(\cA)$.
\end{defn}
Then, one can define an isomorphism between those two algebras on the level of universal elements, which descends directly to an isomorphism on the elements.
\begin{prop}\label{ch03:prop-drinfeld-loop-isomorphism}
	There exists an isomorphism $\psi: \bar{\mathcal{L}} \to D(\cA)$ defined by the equations
	\begin{equation}
	(id\otimes \psi) M_1= N_+ G , \qquad (id\otimes \psi) M_0 = N_+ N_-^{-1} ,
	\end{equation}	
	between the Drinfeld double $D(\cA)$ and the algebra $\bar{\mathcal{L}}$.	
\end{prop}

We provide the proof of this assertion for the graded Hopf algebras below. In order to do that, we first give some additional propositions and lemmas.

\begin{prop}\label{appendixA:prop7}
 The universal element $G \in \cA\otimes \cA^*\subset \cA \otimes D(\cA)$ defined in Definition \ref{ch03:def-universal-element-G} satisfies the equation
 \begin{equation}
  (\Delta \otimes id) G = G_{13} G_{23} .
 \end{equation}
\begin{proof}
 Analogous to the proof of Proposition \ref{appendixA:prop1}.
\end{proof}
\end{prop}

\begin{prop}\label{appendixA:prop8}
 The crossing relations \eqref{susytt} are equivalent to the requirement that for all $\xi\in\cA$ the universal element $G \in \cA\otimes \cA^* \subset \cA \otimes D(\cA)$ ought to satisfy
 \begin{equation}
  G (id\otimes \iota)\Delta(\xi) = (id\otimes \iota)\Delta^{op}(\xi) G,
 \end{equation}
 where $\iota: \cA \to D(\cA)$ is a canonical embedding of $\cA$ in $D(\cA)$.
\begin{proof}
 Setting $\xi = e_\alpha$ and dropping $\iota$ for clarity one can calculate using the decomposition of the element $G$
 \begin{align*}
  \text{LHS} &= \sum_{\rho\in I} (-1)^{|\rho|} \Delta^{op}(e_\alpha) (e_\rho \otimes e^\rho) = \sum_{\gamma,\rho,\sigma\in I} (-1)^{|\gamma||\sigma|+|\rho|} \mu^{ \sigma \gamma}_{\alpha} (e_\gamma \otimes  e_\sigma) (e_\rho\otimes  e^\rho) = \\
  &= \sum_{\gamma,\rho,\sigma\in I} (-1)^{|\gamma||\sigma|+|\rho||\sigma|+|\rho|} \mu^{ \sigma \gamma}_{\alpha}e_\gamma e_\rho \otimes e_\sigma  e^\rho = = \sum_{\beta,\gamma,\rho,\sigma\in I} (-1)^{|\gamma||\sigma|+|\rho||\sigma|+|\rho|} \mu^{ \sigma \gamma}_{\alpha} m_{\gamma \rho}^{\beta} e_\beta \otimes e_\sigma  e^\rho = \\
  &= \sum_{\beta\in I} e_\beta \otimes \left[ \sum_{\gamma,\rho,\sigma\in I} (-1)^{|\beta||\sigma|+|\rho|} \mu^{ \sigma \gamma}_{\alpha} m_{\gamma \rho}^{\beta} e_\sigma  e^\rho \right] ,\\
  \text{RHS} &= \sum_{\rho\in I} (-1)^{|\rho|} (e_\rho \otimes e^\rho) \Delta(e_\alpha) = \sum_{\gamma,\rho,\sigma\in I} (-1)^{|\rho|} \mu^{\gamma \sigma}_{\alpha} (e_\rho \otimes e^\rho)(e_\gamma \otimes e_\sigma) = \\
  & = \sum_{\gamma,\rho,\sigma\in I} (-1)^{|\gamma||\rho|+|\rho|} \mu^{\gamma \sigma}_{\alpha}  e_\rho e_\gamma \otimes e^{\rho}e_{\sigma} = \sum_{\beta,\gamma,\rho,\sigma\in I} (-1)^{|\gamma||\rho|+|\rho|} m_{\rho \gamma}^{\beta} \mu^{\gamma \sigma}_{\alpha} e_\beta \otimes e^{\rho}e_{\sigma}=\\
  & = \sum_{\beta\in I} e_\beta \otimes \left[ \sum_{\gamma,\rho,\sigma\in I} (-1)^{|\rho||\gamma|+|\rho|} m_{\rho \gamma}^{\beta} \mu^{\gamma \sigma}_{\alpha} e^{\rho}e_{\sigma} \right] .
 \end{align*}
 Because the basis elements $e_\beta$ are linearly independent, the elements in the second tensor factor have to be equal for each $\beta$ independently, which leads to the relation
 \begin{equation*}
  \sum_{\gamma,\rho,\sigma\in I} (-1)^{|\beta||\sigma|+|\rho|} \mu^{\sigma\gamma}_{\alpha} m^\beta_{\gamma\rho} e_\sigma e^\rho = \sum_{\gamma,\rho,\sigma\in I} (-1)^{|\rho||\gamma| + |\rho|} m^\beta_{\rho\gamma} \mu^{\gamma\sigma}_\alpha e^\rho e_\sigma,
 \end{equation*}
 which was to be proven.
\end{proof}
\end{prop}

\begin{lem}\label{appendixA:lemma-quasitriangular-N-drinfeld}
 The universal elements $N_+$, $N_-$ given in Definition \ref{ch03:def-universal-element-N} satisfy
 \begin{equation*}
  (\Delta\otimes id) (N_\pm) = (N_\pm)_{23} (N_\pm)_{13} .
 \end{equation*}
 \begin{proof}
  Analogous to the proof of Lemma \ref{appendixA:lemma-quasitriangular-N}.
 \end{proof}
\end{lem}

\begin{lem}
 The universal elements $N_+$, $N_-$ satisfy
 \begin{equation*}
  R_{12} (N_\pm)_{23} (N_\pm)_{13} = (N_\pm)_{13} (N_\pm)_{23} R_{12}, \qquad R_{12} (N_+)_{23} (N_-)_{13} = (N_-)_{13} (N_+)_{23} R_{12} .
 \end{equation*}
 \begin{proof}
  Analogous to the proof of Lemma \ref{appendixA:lemma-quasitriangular-N-2}.
 \end{proof}
\end{lem}

\begin{lem}\label{appendixA:quadratic-eq-for-N-lemma-drinfeld}
 The universal element $N = N_+ N_-^{-1}$ satisfies the following relation
 \begin{equation*}
  R(\Delta\otimes id) N = N_{13} R_{12} N_{23} .
 \end{equation*}
 \begin{proof}
  Analogous to the proof of Lemma \ref{appendixA:quadratic-eq-for-N-lemma}.
 \end{proof}
\end{lem}

\begin{prop}\label{appendixA:prop12}
 The universal elements $N_+$, $N_-$ and $G$ satisfy the following quadratic equations
 \begin{equation}
  (N_+)_{13} R'_{12} G_{23} = G_{23} R'_{12} (N_+)_{13}, \quad (N_-)_{13} R^{-1}_{12} G_{23} = G_{23} R^{-1}_{12} (N_-)_{13},
 \end{equation}
 \begin{proof}
  Let us assume that the $R$ matrix decomposes into elements $f_i, g_i \in \cA$ as follows
  \begin{equation*}
   R = \sum_{i\in I} f_i \otimes g_i, \qquad R' = \sum_{i\in I} (-1)^{|f_i|} g_i \otimes f_i .
  \end{equation*}
  Then, using the relation
  \begin{equation*}
   (id \otimes \Delta^{op})(R') = R'_{13} R'_{12} ,
  \end{equation*}
  which follows from the left equation \eqref{ch3:eq-quasi-triangularity}, as well as using Proposition \ref{appendixA:prop8} one can prove that
  \begin{align*}
   G_{23} R'_{12} (N_+)_{13} &= G_{23} \Big[(id\otimes id\otimes\iota) R'_{12} R'_{13} \Big] = G_{23} \Big[(id\otimes id\otimes\iota) (id \otimes \Delta)(R') \Big] = \\
   &= G_{23} \sum_{i\in I} (-1)^{|f_i|} (g_i \otimes (id\otimes\iota) \Delta(f_i))  = \\
   &= \sum_{i\in I} (-1)^{|f_i|} (g_i \otimes (id\otimes\iota) \Delta^{op}(f_i)) G_{23} = \\
   &= \Big[(id\otimes id\otimes\iota) (id \otimes \Delta^{op})(R') \Big] G_{23} = \Big[(id\otimes id\otimes\iota) R'_{13} R'_{12} \Big] G_{23} = \\
   &= (N_+)_{13} R'_{12} G_{23} .
  \end{align*}
 The second equation, with the use of the identity
  \begin{equation*}
   (id \otimes \Delta^{op})(R^{-1}) = R^{-1}_{13} R^{-1}_{12} ,
  \end{equation*}
 which follows from the right equation \eqref{ch3:eq-quasi-triangularity} can be show analogously.
 \end{proof}
\end{prop}

At last, we can prove Proposition \ref{ch03:prop-drinfeld-loop-isomorphism}.

\begin{proof}[Proof (of Proposition \ref{ch03:prop-drinfeld-loop-isomorphism})]

In order to do that, we need to show that the isomorphism $\psi$ preserves the defining relations \eqref{ch03:loop-algebra-eq-1}-\eqref{ch03:loop-algebra-eq-3}.  The equation \eqref{ch03:loop-algebra-eq-1} is a direct consequence of Lemma \ref{appendixA:quadratic-eq-for-N-lemma-drinfeld}.

For the equation \eqref{ch03:loop-algebra-eq-2} we can perform the following chain of application of identities given above
 \begin{align*}
  R_{12} (\Delta \otimes \psi) M_1 &= R_{12} (\Delta \otimes id)(N_+ G) = R_{12} \Big[ \underbrace{(\Delta \otimes id)N_+}_{=(N_+)_{23} (N_+)_{13}}\Big] \Big[ \underbrace{(\Delta \otimes id)G}_{=G_{13}G_{23}}\Big] = \\
  &= \underbrace{R_{12} (N_+)_{23} (N_+)_{13}}_{=(N_+)_{13} (N_+)_{23} R_{12}} G_{13}G_{23} = (N_+)_{13} \underbrace{(N_+)_{23} R_{12} G_{13}}_{=G_{13}R_{12}(N_+)_{23}} G_{23} = \\
  &= (N_+)_{13} G_{13}R_{12}(N_+)_{23} G_{23} = (id\otimes id\otimes \psi) (M_1)_{13} R_{12} (M_1)_{23} .
 \end{align*}
Finally, the equation \eqref{ch03:loop-algebra-eq-3} is shown using the following chain of identities
 \begin{align*}
  &(id\otimes id\otimes \psi) R^{-1}_{12} (M_1)_{13} R_{12} (M_0)_{23} = R^{-1}_{12} (N_+)_{13} \underbrace{G_{13} R_{12} (N_+)_{23}}_{=(N_+)_{23}R_{12}G_{13}} (N_-)^{-1}_{23} = \\
  &= \underbrace{R^{-1}_{12} (N_+)_{13} (N_+)_{23}}_{=(N_+)_{23} (N_+)_{13}R^{-1}_{12}} R_{12} \underbrace{G_{13} (N_-)^{-1}_{23}}_{=R'_{12} (N_-)^{-1}_{23} G_{13} (R')^{-1}_{12}} = (N_+)_{23} \underbrace{(N_+)_{13} R'_{12} (N_-)^{-1}_{23}}_{=(N_-)^{-1}_{23}R'_{12}(N_+)_{13}} G_{13} (R')^{-1}_{12} =\\
  &= (N_+)_{23} (N_-)^{-1}_{23} R'_{12} (N_+)_{13} G_{13} (R')^{-1}_{12} =  (id\otimes id\otimes \psi) (M_0)_{23} R'_{12} (M_1)_{23} (R')^{-1}_{12} .
 \end{align*}
 This completes the proof of the proposition.
 \end{proof}

\section{Heisenberg and Drinfeld doubles of the Borel half of $U_q(gl(1|1))$ ($q$ a root of unity)}\label{sec:borel-half-gl11-heisenberg-drinfeld}

In this section, we illustrate general constructions of Heisenberg and Drinfeld doubles when applied to the Borel half of $U_q(gl(1|1))$ when the quantisation parameter $q$ is a root of unity. We first consider the Heisenberg double, and then proceed to the Drinfeld double.

\subsection{Heisenberg double of the Borel half of $U_q(gl(1|1))$} \label{chapter2.3}

We will consider a Heisenberg double of the Borel half of restricted $U_q(gl(1|1))$. The Borel half $\cA = \cB(U_q(gl(1|1))$ of $U_q(gl(1|1))$ is generated by two even elements $\ka, \kb$ and one odd element $e_+$ which are subject to the (anti-)commutation relations
\begin{equation}\label{eq:borel-half-gl11-commutation-relations}
 \begin{gathered}
 \ka^p=\kb^p=1, \qquad \ka \kb = \kb \ka,\\
 \ka e_+ = e_+ \ka, \qquad \kb e_+ = q e_+ \kb, \qquad \{e_+,e_+\} = 0 ,
 \end{gathered}
\end{equation}
where $q=e^{2\pi i/p}$ is a root of unity and $p\in2\mathbb{Z}_{>0}+1$ is an odd positive integer, a co-product
\begin{align}\label{eq:borel-half-gl11-coproduct}
\Delta(\ka) = \ka\otimes\ka, \quad \Delta(\kb) = \kb\otimes\kb, \qquad \Delta(e_+)= e_+\otimes1 + \ka^{-1}\otimes e_+,
\end{align}
and the antipode
\begin{align}\label{eq:borel-half-gl11-antipode} 
 \gamma(\ka) = \ka^{-1}, \qquad \gamma(\kb) &= \kb^{-1}, \qquad \gamma(e_+) = -\ka e_+. 
\end{align}
Because of the first and last relation, the algebra $\cA$ is finite-dimensional --- it has $2p^2$ linearly-independent elements and is spanned by a basis
\begin{equation}\label{glbasis}
  e_{n,m,r} = \ka^n \kb^m e_+^r ,
\end{equation}
where $n,m=0,\ldots,p-1$ and $r=0,1$.

With the use of multiplication and co-multiplication rules for the generators, one can obtain the multiplication and co-multiplication coefficients for the basis elements $e_{n,m,r}$ which are as follows
\begin{align*}
 m_{n_1,m_1,r_1,n_2,m_2,r_2}^{n,m,r} &= q^{-m_2 r_1} \delta_{n,n_1+n_2} \delta_{m,m_1+m_2} \delta_{r,r_1+r_2} \Theta(1-r) \Theta(1-r_1) \Theta(1-r_2), \\
 \mu^{n_1,m_1,r_1,n_2,m_2,r_2}_{n,m,r} &= \delta_{n-r_2,n_1}\delta_{n,n_2} \delta_{m,m_1}\delta_{m,m_2} \delta_{r,r_1+r_2} \Theta(1-r) \Theta(1-r_1) \Theta(1-r_2).
\end{align*}

Now one can move to the dual algebra $\cA^*$. It is generated by two even elements $\hat ka, \hkb$ and one odd element $e_-$ which are subject to the (anti-)commutation relations
\begin{equation}\label{eq:borel-half-gl11-dual-commutation-relations}
 \begin{gathered}
 \hka^p=\hkb^p=1, \qquad \hka \hkb = \hkb \hka, \\ 
 \hka e_- = e_- \hka, \qquad \hkb e_- = q^{-1} e_- \hkb, \qquad \{e_-,e_-\} = 0 ,
 \end{gathered}
\end{equation}
where $q=e^{2\pi i/p}$ is the same root of unity as before, co-product
\begin{align}\label{eq:borel-half-gl11-dual-coproduct}
\hat \Delta(\hka) = \hka\otimes\hka, \quad \hat \Delta(\hkb) = \hkb\otimes\hkb, \qquad \hat \Delta(e_-)= e_-\otimes1 + \hka\otimes e_-,
\end{align}
and the antipode is
\begin{align} \label{eq:borel-half-gl11-dual-antipode}
 \hat \gamma(\hka) = \hka^{-1}, \qquad \hat \gamma(\hkb) &= \hkb^{-1}, \qquad \hat \gamma(e_-) = - e_- \hka^{-1}. \end{align}
The algebra $\cA^*$ has $2p^2$ linearly-independent elements and is spanned by a dual basis
\begin{equation}\label{gldualbasis}
  e^{n,m,r} = \frac{1}{p^2} (q-q^{-1})^r \hka^{-r} e_-^r \sum_{t,t'=0}^{p-1} q^{nt'+mt} \hka^{-t} \hkb^{-t'}  ,
\end{equation}
where $n,m=0,\ldots,p-1$ and $r=0,1$.  One can show that the elements $e_{n,m,r}$ and $e^{n,m,r}$ are dual to each other w.r.t. a Hopf pairing
\begin{align}\label{eq:borel-half-gl11-dual-duality-bracket}
 (e_{n_1,m_1,r_1},e^{n_2,m_2,r_2}) = \delta_{n_1,n_2}\delta_{m_1,m_2}\delta_{r_1,r_2}.
\end{align}
The Heisenberg double $H(\cA)$ is then given by the those elements which satisfy in addition the (anti-)commutation relations \eqref{heisenbergmultiplication}
\begin{equation}
\begin{aligned}
 e_{n,m,0} e^{n',m',0} &= e^{n'-n,m'-m,0} e_{n,m,0} ,\\
 e_{n,m,1} e^{n',m',0} &= e^{n'-n+1,m'-m,0} e_{n,m,1} ,\\
 e_{n,m,0} e^{n',m',1} &= q^{-m} e^{n'-n,m'-m,1} e_{n,m,0} ,\\
 e_{n,m,1} e^{n',m',1} &= e^{n'-n,m'-m,0} e_{n,m,0} - q^{-m} e^{n'-n+1,m'-m,1} e_{n,m,1} ,
\end{aligned} 
\end{equation}
where we dropped the trivial tensor products for simplicity of notation. Those relations can be translated into the relations between the generators $\ka,\kb,\hka,\hkb, e_+, e_-$
\begin{equation}\label{eq:commutationrelation_heis_borel_gl11}
\begin{aligned}
 \ka \hka &= \hka \ka, &&  \kb \hkb = \hkb \kb, && \ka \hkb = q \hkb \ka, \\
 \kb \hka &= q \hka \kb, && \ka e_- = e_- \ka, && \kb e_- = q^{-1} e_- \kb, \\
 e_+ \hka &= \hka e_+, && e_+ \hkb = q^{-1} \hkb, e_+ , && \{ e_+, e_-\}= \frac{\hka}{q-q^{-1}}.
\end{aligned} 
\end{equation}
We note that the Heisenberg double $H(\cA)$ has $4 p^4$ elements spanned by elements
\begin{equation}
\label{eq:heis_borel_gl11-basis}
 \{ \ka^n \kb^m \hka^{n'} \hkb^{m'} e_+^r e_-^s \}_{n,n',m,m'=0; r,s=0}^{p-1;1} .
\end{equation}
The canonical element $W$ is given in terms of generators as
\begin{align}
 &W = \left(1 - (q-q^{-1}) e_+\otimes e_- \right) \frac{1}{p^2} \sum_{n,m,n',m'=0}^{p-1} q^{nm'+mn'} \ka^n \kb^m \otimes \hka^{-n'} \hkb^{-m'} .
\end{align}

\subsection{Drinfeld double of the Borel half of $U_q(gl(1|1))$} \label{chapter3.3}

We will now consider a Drinfeld double of the Borel half of restricted $U_q(gl(1|1))$. The Borel half $\cA = \cB(U_q(gl(1|1))$ of $U_q(gl(1|1))$ is generated by two even elements $\ka, \kb$ and one odd element $e_+$ which are subject to the relations \eqref{eq:borel-half-gl11-commutation-relations}-\eqref{eq:borel-half-gl11-antipode}, with basis elements $e_{n,m,r}$ given by \eqref{glbasis}.

The dual algebra $\cA^*$ is generated by two even elements $\hat ka, \hkb$ and one odd element $e_-$ which are subject to the relations \eqref{eq:borel-half-gl11-dual-commutation-relations}-\eqref{eq:borel-half-gl11-dual-antipode}, with basis elements $e^{n,m,r}$ given by \eqref{gldualbasis}. The elements $e_{n,m,r}$ and $e^{n,m,r}$ are dual to each other w.r.t. a Hopf pairing \eqref{eq:borel-half-gl11-dual-duality-bracket}.

The Drinfeld double $D(\cA)$ is then given by the those elements which satisfy the (anti-)commutation relations \eqref{drinfeld_definition_eq1}
\begin{equation}
\begin{aligned}
 e_{n,m,0}e^{r,s,0} &=  e^{r,s,0}e_{n,m,0}, \\
 e_{n,m,1} e^{r,s,0} &= e^{r+1,s,0}e_{n,m,1} \\
 e_{n,m,0} e^{r,s,1} &= q^{-m} e^{r,s,1}e_{n,m,0}, \\
 e_{n,m,1} e_{r,s,1} &= -q^{-m}  e^{r+1,s,1}e_{n,m,1} + e^{r,s,0}e_{n,m,0} - q^{-s} e_{n-1,m,0} e^{r,s,0},  
\end{aligned} 
\end{equation}
where we dropped the trivial tensor products for simplicity of notation. Those relations can be translated into the relations between the generators $\ka,\kb,\hka,\hkb, e_+, e_-$
\begin{equation}
\begin{aligned}
\ka \hka &= \hka \ka, &&  \kb \hkb = \hkb \kb, && \ka \hkb = \hkb \ka, \\
\kb \hka &= \hka \kb, && \ka e_- = e_- \ka, && \kb e_- = q^{-1} e_- \kb, \\
e_+ \hka &= \hka e_+, && e_+ \hkb = q^{-1} \hkb, e_+ , && \{ e_+, e_-\}= \frac{\hka - \ka^{-1}}{q-q^{-1}},
\end{aligned} 
\end{equation}
with co-product
\begin{align}
\begin{aligned}
\Delta(\ka) &= \ka\otimes\ka, && \Delta(\kb) = \kb\otimes\kb, && \Delta(e_+)= e_+\otimes1 + \ka^{-1}\otimes e_+, \\
\Delta(\hka) &= \hka\otimes\hka, && \Delta(\hkb) = \hkb\otimes\hkb, && \Delta(e_-)= e_-\otimes\hka + 1\otimes e_-,
\end{aligned} 
\end{align}
and the antipode
\begin{align} 
\begin{aligned}
\gamma(\ka) &= \ka^{-1}, && \gamma(\kb) = \kb^{-1}, && \gamma(e_+) = -\ka e_+ , \\
\gamma(\hka) &= \hka^{-1}, && \gamma(\hkb) = \hkb^{-1}, && \gamma(e_-) = - e_- \hka^{-1}. 
\end{aligned} 
\end{align}
The Drinfeld double $D(\cA)$ has $4 p^4$ elements spanned by elements
\begin{equation}
\{ \ka^n \kb^m \hka^{n'} \hkb^{m'} e_+^r e_-^s \}_{n,n',m,m'=0; r,s=0}^{p-1;1} .
\end{equation}
The universal $R$-matrix is given by
\begin{align}
 &R = \left(1 - (q-q^{-1}) e_+\otimes e_- \right) \frac{1}{p^2} \sum_{n,m,n',m'=0}^{p-1} q^{nm'+mn'} \ka^n \kb^m \otimes \hka^{-n'} \hkb^{-m'} . 
\end{align}
From the corollary \ref{sec3:corollary:drinfeld-double-algebra-homomorphism} one can derive an algebra homomorphism $\phi: D(\cA) \to H(\cA)\otimes H(\cA^{op,cop})$, which for generators of $U_q(gl(1|1))$ take the form
\begin{equation}
 \begin{aligned}
  \phi^{(a,b)}(\ka) &= \ka\otimes\ka^{-1}, \\
  \phi^{(a,b)}(\kb) &= \kb\otimes\kb^{-1}, \\
  \phi^{(a,b)}(e_+) &= (-1)^a e_+\otimes 1 + (-1)^{b+1} \ka^{-1}\otimes \ka e_+, \\
  \phi^{(a,b)}(\hka) &= \hka\otimes\hka^{-1}, \\
  \phi^{(a,b)}(\hkb) &= \hkb\otimes\hkb^{-1}, \\
  \phi^{(a,b)}(e_-) &= (-1)^a e_-\otimes\hka^{-1} + (-1)^{b+1} 1 \otimes e_- .
 \end{aligned}
\end{equation}
By taking a quotient
\begin{align}
 \hka = \ka, \qquad \hkb = \kb ,
\end{align}
one obtains from the Drinfeld double the standard relations for $U_q(gl(1|1))$, i.e.
\begin{equation}
	U_q(gl(1|1)) = D(\mathcal{B}(U_q(gl(1|1))))/\sim .
\end{equation}

\section{Heisenberg double of $U_q(gl(1|1))$ ($q$ a root of unity)} \label{sec:gl11-heisenberg}
 \label{exmp4}
In this section, we will consider a Heisenberg double of the restricted $U_q(gl(1|1))$. when $q$ is a root of unity. The $\cA=U_q(gl(1|1))$ is generated by two even elements $\ka, \kb$ and two odd element $e_+,e_-$ which are subject to the (anti-)commutation relations
\begin{equation}\label{ch02:eq-gl11-commutation-relations}
\begin{aligned}
 &\ka^p=\kb^p=1, && \ka \kb = \kb \ka,\\
 &\ka e_{\pm} = e_{\pm} \ka, && \kb e_{\pm} = q^{\pm} e_{\pm} \kb, \\
 &\{e_{\pm},e_{\pm}\} = 0, && \{e_{+},e_{-}\} =\frac{\ka-\ka^{-1}}{q-q^{-1}}.
\end{aligned}
\end{equation}
where $q=e^{2\pi i/p}$ is a root of unity and $p\in2\mathbb{Z}_{>0}+1$ is an odd positive integer, and a co-product
\begin{align}\label{ch02:eq-gl11-coproduct}
\begin{aligned}
\Delta(\ka) = \ka\otimes\ka, & && \Delta(\kb) = \kb\otimes\kb, \\
\Delta(e_+)= e_+\otimes1 + \ka^{-1}\otimes e_+, & && \Delta(e_-)= e_-\otimes \ka + 1\otimes e_-,
\end{aligned}
\end{align}
and the antipode is
\begin{align} \gamma(\ka) = \ka^{-1}, \qquad \gamma(\kb) &= \kb^{-1}, \qquad \gamma(e_\pm) = -\ka^{\pm} e_{\pm}. \end{align}
Because of the first and last relation, the algebra $\cA$ is finite-dimensional --- it has $4p^2$ linearly-independent elements and is spanned by a basis
\begin{equation}
e_{n,m,r,s} = \ka^n \kb^m e_+^r e_-^s ,
\end{equation}
where $n,m=0,\ldots,p-1$ and $r,s=0,1$.

With the use of multiplication and co-multiplication rules for the generators, one can obtain the multiplication and co-multiplication coefficients for the basis elements $e_{n,m,r,s}$ which are as follows
\begin{align*}
 m_{n_1,m_1,r_1,s_1,n_2,m_2,r_2,s_2}^{n,m,r,s} &= q^{m_2 (s_1-r_1)} \Theta(1-r) \Theta(1-r_1) \Theta(1-s) \Theta(1-s_2) \times \\
 &\times \Big[ (-1)^{s_1r_2} \delta_{n,n_1+n_2} \delta_{r,r_1+r_2} \delta_{s,s_1+r_2} \Theta(1-r_2) \Theta(1-s_1) +\\
 &+ \frac{1}{q-q^{-1}} (\delta_{s_1,1}\delta_{r_2, 1} \delta_{r,r_1} \delta_{s,s_2}
(\delta_{n,n_1+n_2+1} -\delta_{n,n_1+n_2-1} ) \Big] \delta_{m,m_1+m_2}, \\
\mu^{n_1,m_1,r_1,s_1,n_2,m_2,r_2,s_2}_{n,m,r,s} &= (-1)^{s_1 r_2} \delta_{n_1+r_2,n}\delta_{n_2-s_1,n} \delta_{m,m_1} \delta_{m,m_2} \delta_{r,r_1+r_2} \delta_{s,s_1+r_2} \times \\
 &\times \Theta(1-r) \Theta(1-r_1) \Theta(1-r_2)\Theta(1-s) \Theta(1-s_1) \Theta(1-s_2).
\end{align*}
Now one can move to the dual algebra $\cA^*$. It is generated by two even elements $\la$, $\lb$ and two odd element $\zeta_{\pm}$ which are subject to the (anti-)commutation relations
\begin{equation}\label{ch02:eq-gl11-dual-commutation-relations}
\begin{gathered}
\begin{aligned}
 &\la^p = \lb^p = 1, && \la\lb= \lb\la, \\
 &\la \zeta_{\pm} = \zeta_{\pm} \la, && \lb \zeta_{\pm} = q^{-1} \zeta_{\pm} \lb,
\end{aligned}\\
 \{\zeta_{+},\zeta_{-}\} = \{\zeta_{\pm},\zeta_{\pm}\} = 0.
\end{gathered}
\end{equation}
where $q=e^{2\pi i/p}$ is the same root of unity as before, co-product
\begin{align}\label{ch02:eq-gl11-dual-coproduct}
\hat \Delta(\la) = \la\otimes\la, \quad \hat \Delta(\lb) = \lb\otimes\lb + \lb\zeta_-\otimes \zeta_+\la\lb, \qquad \hat \Delta(\zeta_{\pm})= \zeta_\pm\otimes\la^{\mp1}+1\otimes\zeta_\pm,
\end{align}
and the antipode is
\begin{align} 
\hat \gamma(\la) = \la^{-1}, \qquad \hat \gamma(\lb) &= \lb^{-1} - (q-q^{-1}) \zeta_+ \lb^{-1} \zeta_-, \qquad \hat \gamma(\zeta_{\pm}) = -\la^{\pm1} \zeta_\pm. \end{align}
The algebra $\cA^*$ has $4p^2$ linearly-independent elements and is spanned by a dual basis
\begin{equation}
 e^{n,m,r,s} = \sigma_{r+s} \zeta_+^r \frac{1}{p^2} \sum_{t,t'=0}^{p-1} q^{nt'+mt} \la^{-t} \lb^{-t'} \zeta_-^s ,
\end{equation}
where $n,m=0,\ldots,p-1$ and $r,s=0,1$ and $\sigma_0=\sigma_1=-\sigma_2=1$. One can show that the elements $e_{n,m,r,s}$ and $e^{n,m,r,s}$ are dual to each other w.r.t. a Hopf pairing:
\begin{align}\label{eq:borel-gl11-dual-duality-bracket}
(e_{n_1,m_1,r_1,s_1},e^{n_2,m_2,r_2,s_2}) = \delta_{n_1,n_2}\delta_{m_1,m_2}\delta_{r_1,r_2}\delta_{rs_1,s_2}.
\end{align}
The Heisenberg double $H(\cA)$ is then given by the those elements which satisfy in addition the (anti-)commutation relations
\begin{equation*}
\begin{aligned}
 e_{1,0,0,0} e^{n',m',r',s'} &= e^{n'-1,m',r',s'} e_{1,0,0,0} ,\\
 e_{0,1,0,0} e^{n',m',r',s'} &= q^{s'-r'} e^{n',m'-1,r',s'} e_{0,1,0,0} ,\\
 e_{0,0,1,0} e^{n',m',r',s'} &= (-1)^{r'+s'} e^{n'+1,m',r',s'} e_{0,0,1,0} + \\
 &+ \left[ \frac{1}{q-q^{-1}} (-1)^{r'} \delta_{s',0} \left( e^{n'-1,m',r',1} - e^{n'+1,m',r',1} \right) - \delta_{r',1} e^{n',m',0,s'} \right] e_{0,0,0,0} ,\\
 e_{0,0,0,1} e^{n',m',r',s'} &= (-1)^{r'+s'} e^{n',m',r',s'} e_{0,0,0,1} + \delta_{s',1} e^{n',m',r',0} e_{1,0,0,0} ,
\end{aligned} 
\end{equation*}
where we dropped the trivial tensor products for simplicity of notation. Those relations can be translated into the relations between the generators $\ka,\kb,\la,\lb, e_+, e_-, \zeta_+,\zeta_-$
\begin{equation}
\begin{aligned}
\ka \la &= \la \ka, &&  \ka \lb = q\lb \ka, \\
\kb \lb &= \lb \kb, &&  \kb \la = q\la \kb, \\
\end{aligned} 
\end{equation}
and
\begin{equation}
\begin{aligned}
e_-\la &= \la e_-, &&  e_-\lb = \lb e_-, \\
e_+\la &= \la e_+, && e_+\lb = q^{-1} \lb e_+ + \lb \zeta_- , 
\end{aligned} 
\end{equation}
and
\begin{equation}
\begin{aligned}
 \ka \zeta_\pm &= \zeta_\pm \ka, &&  \kb \zeta_\pm = q^{\mp1} \zeta_\pm \kb,
\end{aligned} 
\end{equation}
and
\begin{equation}
\begin{aligned}
 \{e_-, \zeta_-\}&= \ka, && \{e_+, \zeta_+\} = -1 ,\\
 \{e_-, \zeta_+\} &= 0 , &&  \{e_+, \zeta_-\} = 0 .
\end{aligned} 
\end{equation}
We note that the Heisenberg double $H(\cA)$ has $8 p^4$ elements and is spanned by a basis
\begin{equation}
\{ \ka^n \kb^m \la^{n'} \lb^{m'} e_+^r e_-^s \zeta_+^{r'},\zeta_-^{s'}\}_{n,n',m,m'=0; r,s,r',s'=0}^{p-1;1} .
\end{equation}
The canonical element $W$ is given in terms of generators as
\begin{align}
\begin{aligned}
 W &= \frac{1}{p^2} \sum_{n,m,t,t'=0}^{p-1} q^{nt'+mt} \left( \ka^n \kb^m \otimes \la^{-t} \lb^{-t'} \right) \times \\
 &\times \left(1 - \ka e_+ \otimes \zeta_+ - e_- \otimes \zeta_- - \ka e_+e_- \otimes \zeta_+\zeta_- \right)   
\end{aligned}
\end{align}

Finally, let us comment about the connection between the Heisenberg double of $U_q(gl(1|1))$ and the handle algebra associated to the same super Hopf algebra. The elements $\tilde e_{n,m,r,s}$ and $\tilde e^{n,m,r,s}$ of $H(\cA^*)$ satisfy the same multiplication and co-multiplication relations as their ``untilded'' counterparts, and they are generated by the ``tilded'' generators satisfying the relations \eqref{ch02:eq-gl11-commutation-relations}-\eqref{ch02:eq-gl11-coproduct} and \eqref{ch02:eq-gl11-dual-commutation-relations}-\eqref{ch02:eq-gl11-dual-coproduct} respectively. Then, the Heisenberg double $H(\cA^*)$ is given by the elements which satisfy the (anti-)commutation relations
\begin{equation*}
\begin{aligned}
 \tilde e^{n',m',r',s'} \tilde e_{1,0,0,0} &= \tilde e_{1,0,0,0}  \tilde e^{n'-1,m',r',s'} ,\\
 \tilde e^{n',m',r',s'} \tilde e_{0,1,0,0} &= \tilde e_{0,1,0,0}  \tilde e^{n',m'-1,r',s'} ,\\
 \tilde e^{n',m',r',s'} \tilde e_{0,0,1,0} &= (-1)^{r'+s'} \tilde e_{0,0,1,0} \tilde e^{n',m',r',s'} + q^{-m} \delta_{r',1} \tilde e_{-1,0,0,0} \tilde e^{n',m',0,s'} ,\\
 \tilde e^{n',m',r',s'} \tilde e_{0,0,0,1} &= (-1)^{r'+s'} \tilde e_{0,0,0,1} \tilde e^{n'-1,m',r',s'} + q^{m'} \tilde e_{0,0,0,0} \Big[ (-1)^{-r'} \delta_{s',1} \tilde e^{n',m',r',0} + \\
 &+ \frac{1}{q-q^{-1}} \delta_{r',0} \left( \tilde e^{n'-1,m',1,s'} - \tilde e^{n'+1,m',1,s'} \right) \Big] .
\end{aligned} 
\end{equation*}
They translate into the commutation for the generators which have the form
\begin{equation}\label{tildealgebra1}
\begin{aligned}
\tla \tka &= \tka \tla, && \tla \tkb = q \tkb \tla, \\
\tlb \tkb &= \tkb \tlb, && \tlb \tka = q \tka \tlb, \\
\end{aligned} 
\end{equation}
and
\begin{equation}\label{tildealgebra2}
\begin{aligned}
 \tla \tilde e_+ &= \tilde e_+ \tla , && \tlb \tilde e_+ = \tilde e_+ \tlb , \\
 \tla \tilde e_- &= \tilde e_- \tla , && \tlb \tilde e_- = \tilde e_- \tlb + \tilde\zeta_+ \tla \tlb , 
\end{aligned} 
\end{equation}
as well as
\begin{equation}\label{tildealgebra3}
\begin{aligned}
 \tilde\zeta_\pm \tka &= \tka \tilde\zeta_\pm , &&  \tilde\zeta_\pm \tkb = \tkb \tilde\zeta_\pm ,
\end{aligned} 
\end{equation}
and the anti-commutation relations
\begin{equation}\label{tildealgebra4}
\begin{aligned}
 \{\tilde\zeta_-, \tilde e_-\}&= \tla, && \{\tilde\zeta_+, \tilde e_+\} = \tka^{-1} \tla^{-1} ,\\
 \{\tilde\zeta_+, \tilde e_-\} &= 0 , && \{\tilde\zeta_-, \tilde e_+\} = 0 .
\end{aligned} 
\end{equation}

Under the isomorphism $\varphi$ from Proposition \ref{ch02:prop-heisenberg-handle-isomorphism}, the above commutation relations are equivalent to the (anti-)commutation relations for the generators of the handle algebra $\cT(U_q(gl(1|1))$ obtained in \cite{AGPS1}, which are as follows: the handle algebra is generated by the elements $ \ka^{(a)}$, $\kb^{(a)}$, $e_+^{(a)}$, $e_-^{(a)}$, $\ka^{(b)}$, $\kb^{(b)}$, $e_+^{(b)}$, $e_-^{(b)}$, and the generators $\ka^{(a)}$, $\kb^{(a)}$, $e_+^{(a)}$, $e_-^{(a)}$ satisfy the $U_q(gl(1|1))$ (anti-)commutation relations \eqref{ch02:eq-gl11-commutation-relations} between themselves. The same is true of the generators $\ka^{(b)}$, $\kb^{(b)}$, $e_+^{(b)}$, $e_-^{(b)}$. The ``mixed" commutation relations are given by the formulae
\begin{align}\begin{aligned}
&\big(\kb^{(a)}\big)^{2n} \big(\kb^{(b)}\big)^{2m} = \big(\kb^{(b)}\big)^{2m} \big(\kb^{(a)}\big)^{2n} , \\ &\big(\kb^{(a)}\big)^{2n} \big(\ka^{(b)}\big)^{2m} = q^{2nm} \big(\ka^{(b)}\big)^{2m} \big(\kb^{(a)}\big)^{2n} , \\ &\big(\kb^{(a)}\big)^{2n} e_-^{(b)} = q^{-n} e_-^{(b)} \big(\kb^{(a)}\big)^{2n} ,\\ &\big(\kb^{(b)}\big)^{2n} e_+^{(a)} = q^{2n} e_+^{(a)} \big(\kb^{(b)}\big)^{2n} , \\ & \\ &\big[\big(\ka^{(a)}\big)^n, e_+^{(b)}\big] = 0, \\ &\big[\big(\ka^{(a)}\big)^n, e_-^{(b)}\big] = 0, 
\end{aligned} && \begin{aligned}
&\big(\ka^{(a)}\big)^{2n} \big(\ka^{(b)}\big)^{2m} = \big(\ka^{(b)}\big)^{2m} \big(\ka^{(a)}\big)^{2n} , \\ &\big(\ka^{(a)}\big)^{2n} (\kb^{(b)})^{2m} = q^{2nm} (\kb^{(b)})^{2m} (\ka^{(a)})^{2n} , \\ &\big[\big(\kb^{(a)}\big)^{2n}, e_+^{(b)}\big] = q^{n} [n]_q(q-q^{-1}) e_+^{(a)} \big(\kb^{(a)}\big)^{2n}, \\ &\big(\kb^{(b)}\big)^{2n} e_-^{(a)} = q^{-n} e_-^{(a)} \big(\kb^{(b)}\big)^{2n} + \\ &- q^{-2n}[n]_q (q-q^{-1}) (\ka^{(a)}) \big(\ka^{(b)}\big)^{-1} e_-^{(b)} (\kb^{(b)})^{2n} , \\ &\big[\big(\ka^{(b)}\big)^n, e_+^{(a)}\big] = 0, \\ &\big[\big(\ka^{(b)}\big)^n, e_-^{(a)}\big] = 0,
\end{aligned}
\end{align}
and the anti-commutation relations by
\begin{align}\begin{aligned}
& \{e_+^{(a)}, e_+^{(b)} \} = 0,~~~~~~~~~&&\{e_+^{(a)}, e_-^{(b)}\} = \ka^{(b)} (q-q^{-1})^{-1},\\
&\{e_-^{(a)}, e_-^{(b)}\} = 0,
~~~~~~~~~&& \{e_-^{(a)}, e_+^{(b)} \} = \left( \ka^{(a)} - \big(\ka^{(a)}\big)^{-1} - \ka^{(a)}\big (\ka^{(b)}\big)^{-2}\right) (q-q^{-1})^{-1},
\end{aligned}
\end{align}
where the powers $n,m=0,\ldots,p-1$. The generators $\ka^{(a)}$, $\kb^{(a)}$, $e_+^{(a)}$, $e_-^{(a)}$ have a geometrical interpretation as associated to the a-cycle of a torus, while $\ka^{(b)}$, $\kb^{(b)}$, $e_+^{(b)}$, $e_-^{(b)}$ are associated to the b-cycle.

As a final remark, given that the isomorphism $\varphi$ uses the universal element $G$, we would like to explicitly state that in the case of $U_q(gl(1|1))$ this universal element \eqref{universalG} has the following form
\begin{equation}
\begin{aligned}
 G &= \frac{1}{p^2} \sum_{n,m,\tilde n,\tilde m=0}^{p-1} q^{ n\tilde n + m\tilde m} \left( \tka^{n} \tkb^{m} \otimes
 \tlb^{-\tilde n} \tla^{-\tilde m} \right) \times \\
&\times (1\otimes1 + \tka \tilde e_+\otimes\tilde\xi_+ + \tilde e_-\otimes\tilde \xi_- - \tka \tilde e_- \tilde e_+ \otimes \tilde \xi_-\tilde\xi_+ ) .
\end{aligned}
\end{equation}

This universal element $G$ has been first obtained from directly solving the quadratic relation \eqref{appendixA:eq-universal-element-G-relation} in \cite{AGPS2}. Use of the defining equation \eqref{universalG} gives the same result, as one would expect.

\section{Heisenberg and Drinfeld doubles of the Borel half of $U_q(osp(1|2))$ ($q$ a root of unity)}\label{sec:borel-half-osp-heisenberg-drinfeld}

In this section, we illustrate general constructions of Heisenberg and Drinfeld doubles when applied to the Borel half of $U_q(osp(1|2))$ when the quantisation parameter $q$ is a root of unity. We first consider the Heisenberg double, and then proceed to the Drinfeld double.

\subsection{Heisenberg double of the Borel half of $U_q(osp(1|2))$}\label{chapter2.2}

First, lets us consider a finite-dimensional Heisenberg double of the Borel half of $U_q(osp(1|2))$. The Borel half algebra $\mathcal{A}=\mathcal{B}(osp(1|2)))) = \text{span}\{ k^m \vp^n \}_{m,n=0}^{p-1}$ is generated by an even graded element $k$ and an odd graded element $\vp$ with a commutation relation
\begin{align}\label{eq:borel-half-osp-commutation-relations}
 k^p=1, \qquad \vp^p = 0, \qquad k \vp = q \vp k,
\end{align}
a co-product
\begin{align}\label{eq:borel-half-osp-coproduct}
	&\Delta(k) = k\otimes k , &&  \Delta(\vp) = \vp\otimes k^{-1} + 1\otimes \vp, 
\end{align}
where $q = e^{i\pi \ub^2}$ is the deformation parameter which is not a root of unity, the antipode as follows
\begin{align}\label{eq:borel-half-osp-antipode}
&\gamma(k)= k^{-1},~~~~~&& \gamma(\vp) = - \vp k .
\end{align}
We can choose the basis elements of $\mathcal{A}$ in the following way
\begin{equation}\label{eq:borel-half-osp-basis}
	e_{m,n} = k^m \vp^n ,
\end{equation}
and extend the the first index so that it is taken to be modulo $p$ 
\begin{equation}
	e_{m + rp,n} = e_{m,n} .
\end{equation}
Using the properties of the generators $k, \vp$ and the binomial and q-binomial formulae
\begin{align}\label{eq:binomial-q-binomial-formulae}
  (g + h)^n = \sum_{k=0}^n {n\choose k} g^{n-k} h^k, \qquad (u + v)^n = \sum_{k=0}^n {n\choose k}_q u^{n-k} v^k,
\end{align}
where $gh = hg$ and $uv = q^{-1} vu$, one can find the multiplication and co-multiplication of the basis elements
\begin{align}
 & e_{m,n} e_{k,l} = q^{-nk} e_{m+k, n+l} \Theta(p-1-n-l) ,\label{finite_osp_multiplication_basis} \\
 &\Delta(e_{m,n}) = \sum_{k=0}^{n} {n \choose k}_{-q} e_{m,n-k} \otimes e_{m+k-m, k} . \label{finite_osp_comultiplication_basis}  
\end{align}

The Borel half algebra $\mathcal{A}^*=\mathcal{B}(U_q(osp(1|2)))=\text{span}\{\hk^m \vm^n\}_{m,n=0}^{p-1}$ on the other hand is generated by elements $\hk, \vm$ which satisfy a commutation relation
\begin{align}\label{eq:borel-half-osp-dual-commutation-relations}
 \hk^p=1, \qquad \vm^p = 0, \qquad \hk \vm = q^{-1} \vm \hk,
\end{align}
and have a co-product given by
\begin{align}\label{eq:borel-half-osp-dual-coproduct}
	&{\hat \Delta}(\hk) = \hk\otimes\hk,
	&&{\hat \Delta}(\vm) = \vm\otimes \hk + 1\otimes \vm,
\end{align}
with the antipode
\begin{align}\label{eq:borel-half-osp-dual-antipode}
\begin{aligned}
 \hat\gamma(\hk) &= \hk^{-1},~~~~~&& \hat\gamma(\vm) = - \vm \hk. 
 \end{aligned}
\end{align}
The basis elements for $\mathcal{A}^*$ are given by
\begin{equation}\label{eq:borel-half-osp-dual-basis}
 e^{m,n} = (-1)^{\frac{1}{2}n(n+1)} (-q)_n^{-1} \frac{1}{p} \sum_{r=0}^{p-1} q^{mr} \hk^r \vm^n ,
\end{equation}
where q-numbers $(q)_n$ are defined as $(q)_n = (1-q)\ldots(1-q^n)$ and we extend the first index in the following way
\begin{equation}
	e^{m + rp,n} = e^{m,n} .
\end{equation}
Their multiplication and co-multiplication are as follows
\begin{align}
 & e^{m,n} e^{k,l} = (-1)^{nl} {n+l \choose n}_{-q} \delta_{m,k+n} e^{m,n+l} \Theta(p-1-n-l) , \label{finite_osp_multiplication_basis_dual} \\
 &{\hat \Delta}(e^{m,n}) = \sum_{k,r=0}^m (-1)^{k(n-k)} q^{-(n-k)(m-r)} e^{r,n-k} \otimes e^{m-r, k}  . \label{finite_osp_comultiplication_basis_dual} 
\end{align}
One can show that the elements $e_{m,n}$ and $e^{m,n}$ are dual to each other w.r.t. a Hopf pairing:
\begin{align}\label{eq:borel-half-osp-dual-duality-bracket}
(e_{m_1,n_1},e^{m_2,n_2}) = \delta_{m_1,m_2}\delta_{n_1,n_2}.
\end{align}
By inspection the multiplication and co-multiplication coefficients are equal to
\begin{align*}
 & m^{n'',m''}_{n,m;n',m'} = q^{-m n'} \delta_{n'', n + n' \text{ mod } p} \delta_{m'',m+m'} \Theta(p-1-m)\Theta(p-1-m') \Theta(p-1-m'') , \\
 & \mu_{n'',m''}^{n,m;n',m'} = { m'' \choose m}_{-q} \delta_{n,n''} \delta_{n'',n'+m\text{ mod } p}\delta_{m'',m+m'} \Theta(m''-m)\Theta(m''-m') \Theta(p-1-m'')   , 
\end{align*}
and  $\cA$ and $\cA^*$ are dual as Hopf algebras. The Heisenberg double $H(\cA)$ is then given by the those elements which satisfy the (anti-)commutation relations \eqref{heisenbergmultiplication}. One can calculate the (anti-)commutation relations for the lowest basis elements
\begin{equation}
\begin{aligned}
 e_{1,0} e^{n,m} &= q^{-m} e^{n-1,m} e_{1,0} \\
 e_{0,1} e^{n,0} &= e^{n,0} e_{0,1} \\
 e_{0,1} e^{0,1} &= e^{0,0} e_{-1,0} - e^{0,1} e_{0,1},
\end{aligned} 
\end{equation}
where we dropped the trivial tensor products for simplicity of notation. Those relations can be translated into the relations between the even elements $k$ and $\hk$ and the odd elements $\vp$ and $\vm$ satisfying (anti-)commutation relations
\begin{align}\label{finite_osp_heisenberg_commutation_relations}\begin{aligned}
		& k \hk = q^{-1} \hk k, && \{\vp, \vm\} = -(1+q) k^{-1}, \\
		& k\vp = q \vp k, && k\vm = q^{-1} \vm k, \\
		& [\hk, \vp] = 0, && \hk\vm = q^{-1} \vm\hk.
\end{aligned}
\end{align}
The canonical element $W$ is given in terms of generators as
\begin{align}
	&W = \frac{1}{p} \sum_{n,r}^{p-1} q^{nr} (k^n \otimes \hk^r) (\vp\otimes \vm;-q)^{-1}_\infty ,
\end{align}
where the special function $(x;q)_\infty$, known under a name of a compact  \textit{quantum dilogarithm}, is defined as
\begin{align}\label{specialfunction}
 &(x;q)_\infty^{-1} = \prod_{k=0}^\infty \frac{1}{1-xq^k}=  \sum_{k=0}^\infty \frac{x^k}{(q)_k} .
\end{align}
Because $(\vp\otimes \vm)^p = 0$, the sum above has only a finite number of terms. Using the properties of the quantum dilogarithm function one can check explicitly that the pentagon equation is satisfied. In particular, it reduces to the identity
\begin{align}\label{identity}
 (V;-q)_\infty (U;-q)_\infty = (U;-q)_\infty ((1+q)^{-1} [U,V];-q)_\infty (V;-q)_\infty ,
\end{align}
for $U = 1 \otimes \vp \otimes \vm, V= \vp\otimes \vm \otimes 1$, which was shown to be satisfied by the quantum dilogarithm $(x;q)_\infty$  \cite{Ka3}.
The square brackets denote the commutator, and operators $U$ and $V$ satisfy the following algebraic relations
\begin{align}
&W= UV+qVU,&&[U,W] = [V,W] = 0.
\end{align} 

\subsection{Drinfeld double of the Borel half of $U_q(osp(1|2))$}\label{sec:borel-half-osp-drinfeld}

Let us follow the previous construction with a Drinfeld double of the Borel half of $U_q(osp(1|2))$. The Borel half algebra $\mathcal{A}=\mathcal{B}(osp(1|2)))) = \text{span}\{  k^m \vp^n \}_{m,n=0}^{p-1}$ is generated by an even graded element $k$ and an odd graded element $\vp$ satisfying the relations \eqref{eq:borel-half-osp-commutation-relations}-\eqref{eq:borel-half-osp-antipode}, with basis elements $e_{m,n}$ given by \eqref{eq:borel-half-osp-basis}.

The Borel half algebra $\mathcal{A}^*=\mathcal{B}(U_q(osp(1|2)))=\text{span}\{\hk^m \vm^n\}_{m,n=0}^{p-1}$ is generated by elements $\hk, \vm$ which satisfy the relations \eqref{eq:borel-half-osp-dual-commutation-relations}-\eqref{eq:borel-half-osp-dual-antipode}, with basis elements $e^{m,n}$ defined by \eqref{eq:borel-half-osp-dual-basis}.

The multiplication and co-multiplication coefficients for the basis elements $e_{n,m}$ and $e^{n,m}$ were given in Section \ref{chapter2.2}. $\cA$ and $\cA^*$ are dual as Hopf algebras w.r.t. a Hopf pairing \eqref{eq:borel-half-osp-dual-duality-bracket}.

The Drinfeld double $D(\cA)$ is then given by the crossing relations, which for lowest basis elements have the form
\begin{equation}
 \begin{aligned}
  & e_{1,0} e^{n,0} = e^{n,0} e_{1,0}, \\
  & e_{0,1} e^{n,0} = e^{n-1,0} e_{0,1} , \\
  & e_{1,0} e^{-1,1} = q^{-1} e^{-1,1} e_{1,0}, \\
  & e_{0,1} e^{1,1} + e_{0,0} e^{0,0} = e^{0,0} e_{-1,0} - e^{0,1} e_{0,1},
 \end{aligned}
\end{equation}
where we dropped the tensor products for simplicity. For the generators $k,\hk, \vp, \vm$, one has the following (anti-)commutation relations
\begin{align} 
\begin{aligned}
& k \hk = \hk k, && \{\vp, \vm\} = (1+q) (1-k^{-1}), \\
& k\vp = q \vp k, && k\vm = q^{-1} \vm k, \\
& \hk\vp = q \vp\hk, && \hk\vm = q^{-1} \vm\hk. 
\end{aligned}
\end{align}
with the co-product
\begin{equation}
\begin{aligned}
&\Delta(k) = k\otimes k, && \Delta(\vp) = \vp\otimes k^{-1} + 1\otimes \vp, \\
&{ \Delta}(\hk) = \hk\otimes\hk,
&&{ \Delta}(\vm) = \vm\otimes1 + \hk \otimes \vm,
\end{aligned}
\end{equation}
and the antipode
\begin{align}
\begin{aligned}
\gamma(k) &= k^{-1}, && \gamma(\vp) = - \vp k ,\\
\gamma(\hk) &= \hk^{-1}, && \gamma(\vm) = - \hk^{-1} \vm. 
\end{aligned}
\end{align}
The universal $R$-matrix is given by
\begin{align}
 &R = \frac{1}{p} \sum_{n,r}^{p-1} q^{nr} (k^n \otimes \hk^r) (\vp\otimes \vm;-q)^{-1}_\infty .
\end{align}

From the corollary \ref{sec3:corollary:drinfeld-double-algebra-homomorphism} one can derive an algebra homomorphism $\phi: D(\cA) \to H(\cA)\otimes H(\cA^{op,cop})$, which for generators of $U_q(osp(1|2))$ take the form
\begin{equation}
 \begin{aligned}
  \phi^{(a,b)}(k) &= k\otimes k^{-1}, \\
  \phi^{(a,b)}(\hat k) &= \hat k\otimes\hat k^{-1}, \\
  \phi^{(a,b)}(\vp) &= (-1)^a \vp\otimes 1 + (-1)^{b+1} k^{-1} \otimes \vp k, \\
  \phi^{(a,b)}(\vm) &= (-1)^a {\hat k}\vm\otimes1 + (-1)^{b+1} 1 \otimes \vm\hat k.
 \end{aligned}
\end{equation}

\section{Heisenberg and Drinfeld doubles of the Borel half of $U_q(osp(1|2))$ ($q$ not a root of unity)}\label{sec:borel-half-osp-discrete-heisenberg-drinfeld}

After considering the finite-dimensional version of the Borel half of $U_q(osp(1|2))$, we will illustrate general constructions of Heisenberg and Drinfeld doubles when applied to the Borel half of $U_q(osp(1|2))$ when the quantisation parameter $q$ is a not root of unity, i.e. when the algebra has countably many basis elements. We first consider the Heisenberg double, and then proceed to the Drinfeld double. While this construction can be made mathematically precise, we do it mostly formally.

\subsection{Heisenberg double of the Borel half of $U_q(osp(1|2))$}\label{appendix-osp}

Lets us consider a Heisenberg double of the Borel half of $U_q(osp(1|2))$ first. The Borel half algebra $\mathcal{A}=\mathcal{B}(osp(1|2)))) = \text{span}\{ H^m \vp^n \}_{m,n=0}^\infty$ is generated by an even graded element $H$ and an odd graded element $\vp$ with a commutation relation
\begin{align}\label{eq:borel-half-osp-discrete-commutation-relations}
	&[H,\vp] =-i \ub \vp,
\end{align}
a co-product and the antipode as follows
\begin{align}\label{eq:borel-half-osp-discrete-coproduct}
	&\Delta(H) = H\otimes1 + 1\otimes H, &&  \Delta(\vp) = \vp\otimes e^{\pi \ub H} + 1\otimes \vp, \\
&\gamma(H)= -H,~~~~~&& \gamma(\vp) = - \vp e^{-\pi \ub H} ,\label{eq:borel-half-osp-discrete-antipode}
\end{align}
where $q = e^{i\pi \ub^2}$ is the deformation parameter which is not a root of unity. We can choose the basis elements of $\mathcal{A}$ in the following way
\begin{equation}\label{eq:borel-half-osp-discrete-basis}
	e_{m,n} = (-1)^{n(n+1)/2} \frac{q^{n/2}}{m! (-q)_n} (i \ub^{-1} H)^m (i \vp)^n .
\end{equation}
Using the properties of the generators $H, \vp$ and the binomial and q-binomial formulae \eqref{eq:binomial-q-binomial-formulae} one can find the multiplication and co-multiplication of the basis elements
\begin{align}
 & e_{m,n} e_{k,l} =  (-1)^{nl} \sum_{j=0}^k {m+j\choose j} {n+l\choose l}_{-q} \frac{(-n)^{k-j}}{(k-j)!} e_{m+j,n+l},\label{discrete_osp_multiplication_basis} \\
 &\Delta(e_{m,n}) = \sum_{k=0}^m \sum_{l=0}^n \sum_{p=0}^\infty (-1)^{(n-l)l} {k+p\choose k} (n-l)^p (-\pi i \ub^2)^p e_{m-k,n-l}\otimes e_{k+p,l}. \label{discrete_osp_comultiplication_basis}  
\end{align}
The Borel half algebra $\mathcal{A}^*=\mathcal{B}(U_q(osp(1|2)))=\text{span}\{{\hat H}^m \vm^n\}_{m,n=0}^\infty$ on the other hand is generated by elements $\hat H, \vm$ which satisfy a commutation relation
\begin{align}\label{eq:borel-half-osp-discrete-dual-commutation-relations}
	[\hat H,\vm] = +i \ub\vm,
\end{align}
and have a co-product given by
\begin{align}\label{eq:borel-half-osp-discrete-dual-coproduct}
	&{\hat \Delta}(\hat H) = \hat H\otimes1 + 1\otimes \hat H,
	&&{\hat \Delta}(\vm) = \vm\otimes e^{-\pi \ub {\hat H}} + 1\otimes \vm,
\end{align}
with the antipode
\begin{align}\label{eq:borel-half-osp-discrete-dual-antipode}
\begin{aligned}
 \hat\gamma(\hat H) &= - \hat H,~~~~~&& \hat\gamma(\vm) = - \vm e^{\pi \ub \hat H}. 
 \end{aligned}
\end{align}
The basis elements for $\mathcal{A}^*$ are given by
\begin{equation}\label{eq:borel-half-osp-discrete-dual-basis}
 e^{m,n} = (\pi \ub {\hat H})^m (i v^{(-)})^n , 
\end{equation}
and their multiplication and co-multiplication are as follows
\begin{align}
 & e^{m,n} e^{k,l} = \sum_{j=0}^k {k\choose j} (n)^{k-j} (-\pi i \ub^2)^{k-j} e^{m+j,n+l} , \label{discrete_osp_multiplication_basis_dual} \\
 &{\hat \Delta}(e^{m,n}) = \sum_{k=0}^m \sum_{l=0}^n \sum_{p=0}^\infty {m\choose k} {n\choose l}_{-q} \frac{(-n+l)^p}{p!} e^{m-k,n-l}\otimes e^{k+p,l}. \label{discrete_osp_comultiplication_basis_dual} 
\end{align}
One can show that the elements $e_{m,n}$ and $e^{m,n}$ are dual to each other w.r.t. a Hopf pairing:
\begin{align}\label{eq:borel-half-osp-discrete-dual-duality-bracket}
(e_{m_1,n_1},e^{m_2,n_2}) = \delta_{m_1,m_2}\delta_{n_1,n_2}.
\end{align}
under which each degree subalgebra for $\cA$ is matched by the appropriate degree subalgebra from $\cA^*$.
By inspection the multiplication and co-multiplication coefficients are equal to
\begin{align*}
 & m^{r,s}_{m,n;k,l} = (-1)^{nl} {r\choose r-m} {n+l\choose l}_{-q} \frac{(-n)^{k-r+m}}{(k-r+m)!} \Theta(r-m) \Theta(k-r+m) \delta_{s,n+l}, \\
 & \mu_{r,s}^{m,n;k,l} = (-1)^{nl} {k\choose r-m} (n)^{k-r+m} (-\pi i \ub^2)^{k-r+m} \Theta(r-m) \Theta(k-r+m) \delta_{s,n+l} , 
\end{align*}
and  $\cA$ and $\cA^*$ are dual as Hopf algebras. The Heisenberg double $H(\cA)$ is then given by the those elements which satisfy the (anti-)commutation relations \eqref{heisenbergmultiplication}. One can calculate the (anti-)commutation relations for the lowest basis elements
\begin{equation}
\begin{aligned}
 e_{1,0} e^{1,0} &= e^{0,0} e_{0,0} + e^{1,0} e_{1,0} ,\\
 e_{1,0} e^{0,1} &= e^{0,1} e_{1,0} - e^{0,1} e_{0,0} ,\\
 e_{0,1} e^{1,0} &= e^{1,0} e_{0,1},\\
 e_{0,1} e^{0,1} &= -e^{0,1} e_{0,1} + e^{0,0} \sum_{n=0}^\infty (-\pi i \ub^2)^n e_{n,0} ,
\end{aligned} 
\end{equation}
where we dropped the trivial tensor products for simplicity of notation. Those relations can be translated into the relations between the even elements $H$ and $\hat H$ and the odd elements $\vp$ and $\vm$ satisfying (anti-)commutation relations
\begin{align}\label{discrete_osp_heisenberg_commutation_relations}\begin{aligned}
		& [H,\hat H] = \frac{1}{\pi i}, && \{\vp, \vm\} = e^{\pi \ub H}(q^{\frac{1}{2}}+q^{-\frac{1}{2}}) ,\\
		& [H,\vp]=-i \ub\vp, && [H, \vm] = i \ub \vm, \\
		& [\hat H, \vp] = 0, && [\hat H,\vm]=+i \ub \vm.
\end{aligned}\end{align}
The canonical element $W$ is given in terms of generators as
\begin{align}
	&W = \exp(\pi i H\otimes {\hat H}) (-q^\frac{1}{2} \vp\otimes \vm;-q)^{-1}_\infty.
\end{align}
Using the properties of the quantum dilogarithm function one can check explicitly that the pentagon equation  \eqref{pentagonrelation} is satisfied.

\subsection{Drinfeld double of the Borel half of $U_q(osp(1|2))$}\label{3.2}

Let us follow the previous construction with a Drinfeld double of the Borel half of $U_q(osp(1|2))$. The Borel half algebra $\mathcal{A}=\mathcal{B}(osp(1|2)))) = \text{span}\{ H^m \vp^n \}_{m,n=0}^\infty$ is generated by an even graded element $H$ and an odd graded element $\vp$ satisfying the relations \eqref{eq:borel-half-osp-discrete-commutation-relations}-\eqref{eq:borel-half-osp-discrete-antipode} with basis elements $e_{m,n}$ given by \eqref{eq:borel-half-osp-discrete-basis}.

The Borel half algebra $\mathcal{A}^*=\mathcal{B}(U_q(osp(1|2)))=\text{span}\{{\hat H}^m \vm^n\}_{m,n=0}^\infty$ is generated by elements $\hat H, \vm$ subjected to the relations \eqref{eq:borel-half-osp-discrete-dual-commutation-relations}-\eqref{eq:borel-half-osp-discrete-dual-antipode} with basis elements $e^{m,n}$ given by \eqref{eq:borel-half-osp-discrete-dual-basis}. The multiplication and co-multiplication coefficients for the basis elements $e_{n,m}$ and $e^{n,m}$ were given in Section \ref{appendix-osp}. $\cA$ and $\cA^*$ are dual as Hopf algebras  w.r.t. a Hopf pairing \eqref{eq:borel-half-osp-discrete-dual-duality-bracket}.

The Drinfeld double $D(\cA)$ is then given by the crossing relations, which for lowest basis elements have the form
\begin{equation}
 \begin{aligned}
  & e_{1,0} e^{1,0} + e_{0,0} e^{0,0} = e^{1,0} e_{1,0} + e^{0,0} e_{0,0}, \\
  & e_{0,1} e^{1,0} - \pi i \ub^2 e_{0,1} e^{0,0} = e^{1,0} e_{0,1}, \\
  & e_{1,0} e^{0,1} = e^{0,1} e_{1,0} - e^{0,1} e_{0,0}, \\
  & e_{0,1} e^{0,1} + e_{0,0} \sum_{n=0}^\infty \frac{(-1)^n}{n!} e^{n,0} = -e^{0,1} e_{0,1} + e^{0,0} \sum_{n=0}^\infty (-\pi i \ub^2)^n e_{n,0} ,
 \end{aligned}
\end{equation}
where we dropped the tensor products for simplicity. Then, the generators $H,\hat H, \vp, \vm$ satisfy the (anti-)commutation relations as follows
\begin{align} 
\begin{aligned}
&[H,\vp] = -i\ub\vp, &&\qquad [H,\vm] = +i\ub\vm,\\
&[\hat H, \vp] = -i\ub\vp, &&\qquad [\hat H,\vm] = +i\ub\vm,\\
&[H,\hat H] = 0, &&\qquad \{\vp,\vm\} = (q^\frac{1}{2}+q^{-\frac{1}{2}})(e^{\pi \ub H}-e^{-\pi \ub {\hat H}}) .
\end{aligned}
\end{align}
with the co-product
\begin{equation}
\begin{aligned}
&\Delta(H) = H\otimes1 + 1\otimes H, && \Delta(\vp) = \vp\otimes e^{\pi \ub H} + 1\otimes \vp, \\
&{ \Delta}(\hat H) = \hat H\otimes1 + 1\otimes \hat H,
&&{ \Delta}(\vm) = \vm\otimes1 + e^{-\pi \ub {\hat H}}\otimes \vm,
\end{aligned}
\end{equation}
and the antipode
\begin{align}
\begin{aligned}
\gamma(H) &= -H, && \gamma(\vp) = - \vp e^{-\pi \ub H} ,\\
\gamma(\hat H) &= - \hat H, && \gamma(\vm) = - e^{\pi \ub \hat H} \vm. 
\end{aligned}
\end{align}
The universal $R$-matrix is given by
\begin{align}
 &R = \exp(\pi i H\otimes {\hat H}) (-q^\frac{1}{2} \vp\otimes \vm;-q)^{-1}_\infty ,
\end{align}
which is the same form as has been previously obtained \cite {Kulish:1989sv, Kulish:1991}, where the $U_q(osp(1|2))$ algebra has been obtained by directly quantising classical $osp(1|2)$ while imposing that the coproduct ought to be co-associative, as opposed to Drinfeld double construction.

From the corollary \ref{sec3:corollary:drinfeld-double-algebra-homomorphism} one can derive an algebra homomorphism $\phi: D(\cA) \to H(\cA)\otimes H(\cA^{op,cop})$, which for generators of $U_q(osp(1|2))$ take the form
\begin{equation}
 \begin{aligned}
  \phi^{(a,b)}(H) &= H\otimes1 - 1\otimes H, \\
  \phi^{(a,b)}(\hat H) &= \hat H\otimes1 - 1\otimes \hat H, \\
  \phi^{(a,b)}(\vp) &= (-1)^a \vp\otimes e^{\pi\ub H} + (-1)^{b+1} 1\otimes \vp e^{-\pi \ub H}, \\
  \phi^{(a,b)}(\vm) &= (-1)^a \vm\otimes1 + (-1)^{b+1} e^{-\pi\ub \hat H} \otimes e^{\pi \ub \hat H} \vm.
 \end{aligned}
\end{equation}

\section{Outlook} \label{outlook}

 We would like to present three directions as the follow up of this work: one extending a setting where Heisenberg and Drinfeld doubles play an important role to finite-dimensional graded algebras and two concerned with generalising results to non-demumerably infinite graded Hopf algebras.\\

 The first one is in the context of lattice current algebras \cite{Alekseev:1996jz}, which were introduced as a regularisation of the left and right moving degrees of freedom in the WZW models and, in the continuum limit, reproduce WZW current algebra.

 The Drinfeld double plays a role as a part of them. The representation categories of the lattice and the continuum theory coincide. The lattice current algebras and their representation theory was studied for the non-graded, semisimple Hopf algebras. As such, one could investigate them for graded Hopf algebras, especially given the existing work on generalised WZW \cite{Schomerus:2005bf} for the $GL(1|1)$ gauge group. The lattice current algebras for a finite-dimensional, graded, (possibly) non-semisimple Lie algebra $G$ will be studied in a forthcoming work. \\
 
 The second direction for future work would be to extend those Heisenberg and Drinfeld constructions, in which we focused on the finite-dimensional or countably infinite-dimensional Hopf algebras, to the nondenumerably infinite-dimensional Hopf algebras which we believe have the super Hopf algebras generalisations, in the context of non-compact conformal field theories.

 The quantum deformed algebra $U_q(sl(2, \mathbb{R}))$ is known to possess a very interesting self-dual series of infinite-dimensional representations. From the Heisenberg double of the quantum plane Kashaev obtained a class of representations of the $U_q(sl(2, \mathbb{R}))$, as well as its associated $R$-matrix \cite{Ka1}. Those representations were identified to be the class of infinite-dimensional representations $\mathcal{P}_\alpha$ studied \cite{PT1,PT2} in connection with the Liouville field theory, which constitutes a prototypical non-trivial example of the non-compact conformal field theory \cite{T03, Teschner:2003em}. 
 
 Using the means of harmonic analysis, Ponsot and Teschner investigated representations $\mathcal{P}_\alpha$. They have shown that the relation between the fusion category of the conformal field theory and the representation category of the Hopf algebra holds in the case of the Liouville field theory and the $U_q(sl(2, \mathbb{R}))$. Moreover, the consistency of the bootstrap for the Liouville theory, i.e. the fact that the crossing-symmetry equation is satisfied by the three point correlation function, was verified \cite{PT2,PT1,teschner}. The generators $K,E,F$ of the $U_q(sl(2, \mathbb{R}))$ Hopf algebra are realised as positive, self-adjoint operators. The family $\mathcal{P}_\alpha$ is closed under the tensor product and the calculation of the 3j- and 6j-symbols allowed to make a connection with the fusion matrices of the Liouville theory. 
 
 In the context of $\mathcal{N}=1$ supersymmetric Liouville theory, the attempt to find an $U_q(osp(1|2))$ analogue of the representations $\mathcal{P}_\alpha$ resulted obtaining the 6j-symbols which reproduce only the Neveu-Schwarz sector of the theory \cite{Hadasz:2013bwa,Pawelkiewicz:2013wga}. We conjecture, being informed by \cite{APT}, that the proper analogue of the representations $\mathcal{P}_\alpha$, i.e. the one which will encode the entirety of the structure of $\mathcal{N}=1$ supersymmetric Liouville theory, could be obtained using the Heisenberg double and Drinfeld double constructions for nondenumerably infinite Hopf algebras. It would be interesting to generalise the Heisenberg double and Drinfeld double constructions so that they can be applied to the nondenumerably infinite Hopf algebras.
\\

 Finally, a third possible direction of future research could be to apply Heisenberg and Drinfeld double constructions to nondenumerably infinite-dimensional Hopf algebras in the context of Chern-Simons theory in 3-dimensional space and invariants thereof. 
 
 Similarly to the compact Chern–Simons theory, its non-compact version is expected to produce invariants of 3-manifolds from representations of Hopf algebras. However, in contrast to the compact theory, which employed finite-dimensional representations at roots of unity, the non-compact case requires the use of positive, infinite-dimensional representations. The result considering infinite-dimensional representations of  $U_q(sl(2, \mathbb{R}))$ Hopf algebra exists, e.g. the Teichmüller topological quantum field theory \cite{Andersen:2018pnw}, but, at the moment, positive representation theory for the graded Hopf algebras remains unexplored. In particular, it is expected that all known super quantum group representations, and hence all 3-manifold and knot invariants procured from them, can be obtained via analytic continuation of the positive ones.

\acknowledgments
We are very grateful to J\"org Teschner for explanations, suggestions and many helpful discussions and comments on the draft.
We also like to thank Simon Lentner for useful comments about the draft. We thank Rinat Kashaev, Volker Schomerus and Catherine Meusburger for stimulating discussions.

The work of N.A. was supported by funding from the European Union's Horizon 2020 research and innovation programme under the Marie Sklodowska-Curie grant agreement No. 898759 and ReNewQuantum ERC-2018-SyG No. 810573 hosted by Centre for Quantum Mathematics Department of Mathematics and Computer Science (IMADA), Southern Denmark University (SDU), Denmark and by NCCR SwissMAP, funded by the Swiss National Science Foundation at the University of Geneva, Switzerland. For the final part of this project, N.A. acknowledges partial support of the SNF Grant No. 20002-192080  and hosting of University of Zürich. The work of M.K.P. was supported by the European Research Council (advanced grant NuQFT).

 We also thank and acknowledge the AEC centre in France, University of Bern and Max Planck Institute of Mathematics (MPIM) in Bonn for their hospitality during the first stages of this project.

\appendix

\section{Heisenberg and Drinfeld doubles of the Borel half of $U_q(sl(2))$ ($q$ not a root of unity)}  \label{appendix-sl2}

In this appendix we summarise Heisenberg and Drinfeld double constructions when applied to $U_q(sl(2))$ for $q$ not being a root of unity. This is a thoroughly investigated example --- previously studied both in terms of the Heisenberg double \cite{Ka3} as well as the Drinfeld double \cite{DD, Drinfel'd1988} --- however given its ubiquity it is worthwhile to include it to fix conventions as well as for pedagogical reasons as a ``canonical'' illustration of the constructions used.

\subsection{Heisenberg double of the Borel half of $U_q(sl(2))$} \label{appendix-sl2-Heisenberg}

It can be instructive to apply the Heisenberg double construction to the Borel half of $U_q(sl(2))$. The Borel half algebra $\mathcal{A}=\mathcal{B}(U_q(sl(2))) = \text{span}\{ H^m E^n \}_{m,n=0}^\infty$ is generated by elements $H, E$ with a commutation relation
\begin{align}\label{discrete_multiplication}
&[H,E] = -i\ub E,
\end{align}
and a co-product as follows
\begin{align}\label{coproductsltwo}
&\Delta(H) = H\otimes1 + 1\otimes H, && \Delta(E) = E\otimes e^{2\pi \ub H} + 1\otimes E,
\end{align}
where $q = e^{\pi i \ub^2}$ is the deformation parameter and $q$ is not a root of unity. In addition, the antipode is
\begin{align}\label{eq:borel-half-sl2-discrete-antipode}
\begin{aligned}
\gamma(H) &= -H,~~~ && \gamma(E)= - E e^{-2\pi \ub H}.
\end{aligned}
\end{align}
We can choose the basis elements of $\mathcal{A}$ in the following way
\begin{equation}\label{eq:borel-half-sl2-discrete-basis}
e_{m,n} = \frac{q^{n}}{m! (q^2)_n} (i \ub^{-1} H)^m (i E)^n,
\end{equation}
where q-numbers $(q)_n$ are defined as $(q)_n = (1-q)\ldots(1-q^n)$ and $n,m\in\mathbb{N}$.

Using the properties of the generators $H, E$ and the binomial and q-binomial formulae \eqref{eq:binomial-q-binomial-formulae} one can find the multiplication and co-multiplication of the basis elements
\begin{align*}
e_{m,n} e_{k,l} &= \sum_{j=0}^k {m+j\choose j} {n+l\choose l}_{q^2} \frac{(- n)^{k-j}}{(k-j)!} e_{m+j,n+l}, \qquad n>0,\\
e_{m,0} e_{k,l} &= {m+k\choose k} {n+l\choose l}_{q^2} e_{m+k,l}, \\
\Delta(e_{m,n}) &= \sum_{k=0}^m \sum_{l=0}^{n-1} \sum_{p=0}^\infty {k+p\choose k} (n-l)^p (-2\pi i \ub^2)^p e_{m-k,n-l}\otimes e_{k+p,l} + \\
& + \sum_{k=0}^m e_{m-k,0}\otimes e_{k,n} ,
\end{align*}
where ${n\choose k}$ is an ordinary and ${n\choose k}_q = \frac{(q)_{n}}{(q)_{k}(q)_{n-k}}$ is a q-deformed binomial coefficient.

The Borel half algebra $\mathcal{A}^*=\mathcal{B}(U_q(sl(2)))=\text{span}\{{\hat H}^m F^n\}_{m,n=0}^\infty$ is generated by elements $\hat H, F$ which satisfy a commutation relation
\begin{align}\label{discrete_multiplication_dual}
&[\hat H,F] = +i\ub F,
\end{align}
and have a co-product and the antipode given by
\begin{align}\label{discrete_comultiplication_dual}
&{\hat \Delta}(\hat H) = \hat H\otimes1 + 1\otimes \hat H,
&&{\hat \Delta}(F) = F\otimes e^{-2\pi \ub {\hat H}} + 1\otimes F,\\
&{\hat\gamma}(\hat H)= - \hat H,~~~~~&& {\hat\gamma}(F)= - F e^{2\pi \ub \hat H}.
\label{eq:borel-half-sl2-discrete-dual-antipode}
\end{align}
The basis elements for $\mathcal{A}^*$ are given by
\begin{equation}\label{eq:borel-half-sl2-discrete-dual-basis}
e^{m,n} = (2 \pi \ub {\hat H})^m (i F)^n ,
\end{equation}
and their multiplication and co-multiplication is as follows
\begin{align*}
e^{m,n} e^{k,l} &= \sum_{j=0}^k {k\choose j} (n)^{k-j} (-2\pi i \ub^2)^{k-j} e^{m+j,n+l} , \qquad n > 0, \\
e^{m,0} e^{k,l} &= e^{m+k,l} , \\
{\hat \Delta}(e^{m,n}) &= \sum_{k=0}^m \sum_{l=0}^{n-1} \sum_{p=0}^\infty {m\choose k} {n\choose l}_{q^2} \frac{(-n+l)^p}{p!} e^{m-k,n-l}\otimes e^{k+p,l} + \\
& +\sum_{k=0}^m {m\choose k} e^{m-k,0}\otimes e^{k,n} .
\end{align*}
One can show that the elements $e_{m,n}$ and $e^{m,n}$ are dual to each other w.r.t. a Hopf pairing:
\begin{align}\label{eq:borel-half-sl2-discrete-duality-bracket}
(e_{m_1,n_1},e^{m_2,n_2}) = \delta_{m_1,m_2}\delta_{n_1,n_2}.
\end{align}
The multiplication and co-multiplication coefficients are equal to
\begin{align*}
& m^{r,s}_{m,n;k,l} = {r\choose r-m} {n+l\choose l}_{q^2} \frac{(-n)^{k-r+m}}{(k-r+m)!} \Theta(r-m) \Theta(k-r+m) \delta_{s,n+l}, && n>0,\\
& m^{r,s}_{m,0;k,l} = {r\choose r-m} \delta_{r,m+k} \delta_{s,l}, \\
& \mu_{r,s}^{m,n;k,l} = {k\choose r-m} (n)^{k-r+m} (-2\pi i \ub^2)^{k-r+m} \Theta(r-m) \Theta(k-r+m) \delta_{s,n+l} , && n>0, \\
& \mu_{r,s}^{m,0;k,l} = \delta_{r,m+k} \delta_{s,l} ,
\end{align*}
and $\cA$ and $\cA^*$ are dual as Hopf algebras. The Heisenberg double $H(\cA)$ is then given by the those elements which satisfy the (anti-)commutation relations \eqref{heisenbergmultiplication}. One can calculate the (anti-)commutation relations for the lowest basis elements
\begin{equation}
\begin{aligned}
e_{1,0} e^{1,0} &= e^{0,0} e_{0,0} + e^{1,0} e_{1,0} ,\\
e_{1,0} e^{0,1} &= e^{0,1} e_{1,0} - e^{0,1} e_{0,0} ,\\
e_{0,1} e^{1,0} &= e^{1,0} e_{0,1},\\
e_{0,1} e^{0,1} &= e^{0,1} e_{0,1} + e^{0,0} \sum_{n=0}^\infty (-2\pi i \ub^2)^n e_{n,0} ,
\end{aligned} 
\end{equation}
where we dropped the trivial tensor products for simplicity of notation. Those relations can be translated into the relations between the generators $H,\hat H, E, F$, which are as follows
\begin{align}\label{commutationrelation}
\begin{aligned}
&[H,E] = -i\ub E, &&\qquad [H,F] = +i\ub F,\\
&[\hat H, E] = 0, &&\qquad [\hat H,F] = +ibF,\\
&[H,\hat H] = \frac{1}{2\pi i}, &&\qquad [E,F] = (q-q^{-1})e^{2\pi \ub H} .
\end{aligned}
\end{align}
The canonical element $W$ is given in terms of generators as
\begin{align}
&W = \exp(2 \pi i H\otimes\hat H) (- q E\otimes F;q^2)^{-1}_\infty .
\end{align}
Using the properties of the quantum dilogarithm function one can check explicitly that the pentagon equation is satisfied. In particular, it reduces to the identity
\begin{align}\label{identity-sl2}
(V;q^2)_\infty (U;q^2)_\infty = (U;q^2)_\infty ((1-q^2)^{-1} [U,V];q^2)_\infty (V;q^2)_\infty ,
\end{align}
for $U = -q1\otimes E\otimes F,V=-q E\otimes F\otimes 1$, where the operators $U$ and $V$ satisfy the following algebraic relations
\begin{align}
 [U, UV-q^2VU] = [V, UV-q^2VU] = 0.
\end{align} 

\subsection{Drinfeld double of the Borel half of $U_q(sl(2))$} \label{appendix-sl2-Drinfeld}

We will consider now a Drinfeld double of the Borel half of $U_q(sl(2))$. The building blocks that we use are the same as for the Heisenberg double in the Appendix \ref{appendix-sl2-Heisenberg}. Let us remind that the Borel half algebra $\mathcal{A}=\mathcal{B}(U_q(sl(2))) = \text{span}\{ H^m E^n \}_{m,n=0}^\infty$ is generated by elements $H, E$ subjected to the relations \eqref{discrete_multiplication}-\eqref{eq:borel-half-sl2-discrete-antipode} with the basis elements $e_{m,n}$ defined by \eqref{eq:borel-half-sl2-discrete-basis}

The dual Borel half algebra $\mathcal{A}^*=\mathcal{B}(U_q(sl(2)))=\text{span}\{{\hat H}^m F^n\}_{m,n=0}^\infty$ is generated by elements $\hat H, F$ satisfying relations \eqref{discrete_multiplication_dual}-\eqref{eq:borel-half-sl2-discrete-dual-antipode} with the basis elements $e^{m,n}$ defined by \eqref{eq:borel-half-sl2-discrete-dual-basis}.

The multiplication and co-multiplication coefficients for the basis elements $e_{n,m}$ and $e^{n,m}$ were given in Section \ref{appendix-sl2-Heisenberg}. $\cA$ and $\cA^*$ are dual as Hopf algebras w.r.t. a Hopf pairing \eqref{eq:borel-half-sl2-discrete-duality-bracket}.

The Drinfeld double $D(\cA)$ is then given by the crossing relations, whose for lowest basis elements have the form
\begin{equation}
\begin{aligned}
& e_{1,0} e^{1,0} + e_{0,0} e^{0,0} = e^{1,0} e_{1,0} + e^{0,0} e_{0,0}, \\
& e_{0,1} e^{1,0} - 2\pi i \ub^2 e_{0,1} e^{0,0} = e^{1,0} e_{0,1}, \\
& e_{1,0} e^{0,1} = e^{0,1} e_{1,0} - e^{0,1} e_{0,0}, \\
& e_{0,1} e^{0,1} + e_{0,0} \sum_{n=0}^\infty \frac{(-1)^n}{n!} e^{n,0} = e^{0,1} e_{0,1} + e^{0,0} \sum_{n=0}^\infty (-2\pi i \ub^2)^n e_{n,0} ,
\end{aligned}
\end{equation}
where we dropped the tensor products for simplicity. The generators $H,\hat H, E, F$ satisfy thus the commutation relations
\begin{align} 
\begin{aligned}
&[H,E] = -i \ub E, &&\qquad [H,F] = +i\ub F,\\
&[\hat H, E] = -i\ub E, &&\qquad [\hat H,F] = +i \ub F,\\
&[H,\hat H] = 0, &&\qquad [E,F] = (q-q^{-1})(e^{2\pi \ub H}-e^{-2\pi \ub {\hat H}}) ,
\end{aligned}
\end{align}
with the co-product
\begin{equation}
\begin{aligned}
&\Delta(H) = H\otimes1 + 1\otimes H, && \Delta(E) = E\otimes e^{2\pi \ub H} + 1\otimes E, \\
&{ \Delta}(\hat H) = \hat H\otimes1 + 1\otimes \hat H,
&&{ \Delta}(F) = F\otimes 1 + e^{-2\pi \ub {\hat H}}\otimes F,
\end{aligned}
\end{equation}
and the antipode
\begin{align}
\begin{aligned}
\gamma(H) &= -H, && \gamma(E) = - E e^{-2\pi \ub H} ,\\
\gamma(\hat H) &= - \hat H, && \gamma(F) = - e^{2\pi \ub \hat H} F. 
\end{aligned}
\end{align}
The universal $R$-matrix is given by
\begin{align}
&R = \exp(2 \pi i H\otimes\hat H) (- q E\otimes F;q^2)^{-1}_\infty .
\end{align}

From the corollary \ref{sec3:corollary:drinfeld-double-algebra-homomorphism} one can derive an algebra homomorphism $\phi: D(\cA) \to H(\cA)\otimes H(\cA^{op,cop})$, which for generators of $U_q(sl(2))$ take the form
\begin{equation}
\begin{aligned}
\phi^{(a,b)}(H) &= H\otimes1 - 1\otimes H, \\
\phi^{(a,b)}(\hat H) &= \hat H\otimes1 - 1\otimes \hat H, \\
\phi^{(a,b)}(E) &= E\otimes e^{2\pi\ub H} + 1\otimes E e^{-2\pi \ub H}, \\
\phi^{(a,b)}(F) &= F\otimes1 + e^{-2\pi\ub \hat H} \otimes e^{2\pi \ub \hat H} F.
\end{aligned}
\end{equation}


\newpage
\begin{thebibliography}{999}
	\bibliographystyle{JHEP}

	\bibitem{AGPS1}
	N.~Aghaei, A.~M.~Gainutdinov, M.~Pawelkiewicz and V.~Schomerus, \textit{Combinatorial Quantization of $GL(1|1)$ Chern-Simons Theory I: The Torus},
	[arXiv:1811.09123].
	%
	\bibitem{AGPS2}
	N.~Aghaei, A.~M.~Gainutdinov, M.~Pawelkiewicz and V.~Schomerus, \textit{Combinatorial Quantization of $GL(1|1)$ Chern-Simons Theory II : Higher genus}, to appear.
	%
	\bibitem{Drinfeld:1985rx}
	V.~G.~Drinfeld, {\it {Hopf algebras and the quantum Yang-Baxter equation}},
	Sov. Math. Dokl. \textbf{32} (1985), 254-258.
%
	\bibitem{Jimbo:1985zk}
	M.~Jimbo, {\it {A q difference analog of $U(g)$ and the Yang-Baxter equation},}
	Lett. Math. Phys. \textbf{10} (1985), 63-69.
%
	\bibitem{review}
		P. Etingof, V. Ginzburg, N. Guay, D. Hernandez, and A. Savage, {\it {Twenty-five years of representation theory of quantum groups},} Final reports, BIRS, Banff, (2011).
%
		\bibitem{Majidbook}
S.~Majid,{ \it Foundations of Quantum Group Theory}	publisher: Cambridge University Press.Online,
Book ISBN: 9780511613104.
%
\bibitem{Sw}
M.~Sweedler, {\it{Hopf algebras}}, 
Band 44 von Mathematics lecture note series,
Verlag	W.A. Benjamin, 1969, Book ISB:805392556, 9780805392555.
%
\bibitem{Witten:1989rw}
E.~Witten, 
{\it {Gauge Theories, Vertex Models and Quantum Groups}},  {\em
	Nucl.Phys.} \textbf{B330}, (1990), 285.
%
\bibitem{zam}
A.~Zamolodchikov,
 {\it{Factorized S-Matrices in Two Dimensions as the Exact Solutions of Certain Relativistic Quantum Field Theory Models}}, {\em{Annals of Physics}} \textbf{120}, (1979), 253-291.
%
\bibitem{Yang:1967bm}
	C.~N.~Yang,
	{\it Some exact results for the many body problems in one dimension with repulsive delta function interaction},
	Phys.\ Rev.\ Lett.\  \textbf{19}, (1967), 1312.
%
	\bibitem{Baxter:1972hz}
	R.~J.~Baxter,
	{\it Partition function of the eight vertex lattice model},
	Annals Phys.\  \textbf{70}, (1972), 193, [arXiv:0903.3089].
	
	
	
		\bibitem{DD}
	V.~Drinfeld,	{\it  Quantum groups,}  In A. Gleason, editor, Proceedings of the ICM, Rhode Island, (1987), AMS. 798–820.
%
		\bibitem{Majid}
S.~Majid, {\it {Doubles of quasitriangular Hopf algebras}}, Comm. Algebra 19:11,  (1991), 3061-3073, [arXiv:q-alg/9701002].




%
	\bibitem{Majid0}
S.~Majid, {\it {Some remarks on the quantum double, Quantum groups and physics}}, Czechoslovak J. Phys. \textbf{ 44} (1994), 1059–1071, [arXiv:hep-th/9409056].
%



\bibitem{Drinfel'd1988}
V.G. Drinfeld, \textit{Quantum groups}, J. Math. Sci. \textbf{41}, (1988) 898–915. 





\bibitem{Gould}
M. Gould, R. Zhang, AJ. Bracken,
{\it Quantum double construction for graded Hopf algebras},
Bulletin of the Australian Mathematical Society.
1993;47(3):353-375. 


	

\bibitem{Majid:1999}
S. Majid, \textit{Double-bosonization of braided groups and the construction of $U_q(\mathfrak{g})$}, Math. Proc. Camb. Phil. Soc. \textbf{125} (1999), 151—192, [arXiv:q-alg/9511001].


\bibitem{Majid:1990}
S. Majid, \textit{More examples of bicrossproduct and double cross product Hopf algebras}, Israel J. Math. \textbf{72} (1990), 133—148.

\bibitem{Majid:1997}
S. Majid, \textit{New quantum groups by double-bosonisation}, Quantum groups and integrable systems, II (Prague 1996), Czech. J. Phys. \textbf{47} (1) (1997), 79—90, 	[arXiv:q-alg/9610004].

\bibitem{aziz}
R. Aziz and S.Majid,
{\it Co-double bosonisation and dual bases of $cq[SL2]$ and $cq[SL3]$},
J. Algebra, \textbf{518},
2019, 75-118,
[arXiv:math/1703.03456].




	\bibitem{lu}
L. Jiang-Hua,
{\it{On the Drinfeld double and the Heisenberg double of a Hopf algebra}}, Duke Math. J. \textbf{74} {no. 3},  (1994), 763-776.
%




	\bibitem{Ka3}
	R.~M.~Kashaev, {\it The Heisenberg double and the pentagon relation.}
	Algebra in Analiz, 8:4 (1996), 63–74, [arXiv:q-alg/9503005v1].
	%
				\bibitem{Maillet}
	J.~M.~Maillet,
	{\it {On pentagon and Tetrahedron equations}}, Algebra in Analiz, \textbf{6} (1994), 206-214, [arXiv:hep-th/9312037].
	%
	
	\bibitem{Xu:1999eh}
	P.~Xu,
{\it Quantum groupoids},
	Commun. Math. Phys. \textbf{216} (2001), 539-581,
	[arXiv:math/9905192].
	
		\bibitem{Carlip:2004ba}
	S.~Carlip,
	{\it Quantum gravity in 2+1 dimensions: The Case of a closed universe},
	Living Rev. Rel. \textbf{8} (2005), 1,	[arXiv:gr-qc/0409039].
	
	\bibitem{Alekseev:1994pa}
	A.~Y.~Alekseev, H.~Grosse and V.~Schomerus, \textit{Combinatorial quantisation of the Hamiltonian Chern-Simons theory,}
	Comm.\ Math.\ Phys.\  {\bf 172} (1995), 317-358, [arXiv:hep-th/9403066].
	
	
	%
	\bibitem{Alekseev:1994au}
	A.~Y.~Alekseev, H.~Grosse and V.~Schomerus, \textit{Combinatorial quantisation of the Hamiltonian Chern-Simons theory. 2.,}
	Comm.\ Math.\ Phys.\  {\bf 174} (1995), 561-604, [arXiv:hep-th/9408097].
	%
	
	
	
	\bibitem{Alekseev:1995rn}
	A.~Y.~Alekseev and V.~Schomerus, \textit{Representation theory of Chern-Simons observables,}
	Duke Math. J. {\bf 85} (1996), 447-510, [arXiv:q-alg/9503016].
	
	%
	\bibitem{Alekseev:1996ns}
	A.~Y.~Alekseev and V.~Schomerus, \textit{Quantum moduli spaces of flat connections},
	In Goslar 1996, Group theoretical methods in physics:~70, [arXiv:q-alg/9612037].
	

	
	%
	\bibitem{VS}
	V.~Schomerus, 
	\textit{Deformed gauge symmetry in local quantum physics},
	Habilitation 1998, Institute for the theoretical physics at University of Hamburg.
	%
	
	\bibitem{BR95}
	E. Buffenoir and Ph. Roche, 
	{ \it Two dimensional lattice gauge theory based on a quantum group},
	Commun. Math. Phys. 170.3 (1995), 669–698, [arXiv:hep-th/9405126].
	
	%
	\bibitem{BR96}
	E. Buffenoir and Ph. Roche,
	{\it  Link invariants and combinatorial quantization of hamiltonian Chern-Simons theory}, Commun. Math. Phys. 181.2 (1996) 331–365, [arXiv:q-alg/9507001].
	%
	\bibitem{Fock:1998nu}
	V.~V.~Fock and A.~A.~Rosly, \textit{Poisson structure on moduli of flat connections on Riemann surfaces and r matrix,}
	Am.\ Math.\ Soc.\ Transl.\  {\bf 191} (1999), 67-86, [arXiv:math/9802054].
	%
	
	
		%
	\bibitem{Alekseev:1996jz}
	A.~Y.~Alekseev, L.~D.~Faddeev, J.~Frohlich and V.~Schomerus,
	{\it Representation theory of lattice current algebras,}
	Commun. Math. Phys. \textbf{191} (1998), 31-60, 
	[arXiv:q-alg/9604017].
	
	
		\bibitem{Kitaev:1997wr}
	A.~Y.~Kitaev,
	\textit{Fault tolerant quantum computation by anyons},
	Annals Phys. \textbf{303} (2003), 2-30, [arXiv:quant-ph/9707021].
	%
	\bibitem{Meusburger}
	C. Meusburger,
	\textit{Kitaev Lattice Models as a Hopf Algebra Gauge Theory}. Commun. Math. Phys. 353, (2017), 413–468, [arXiv:1607.01144].
	
	
	
	
	
	
	%
	\bibitem{MeusburgerandWise}
	C. Meusburger and K. Wise,
	{\it  Hopf algebra gauge theory on a ribbon graph}
	Reviews in Mathematical Physics \textbf{33}, no. 5, (2021), [arXiv:1512.03966].
	
		\bibitem{book}
	C. Gómez, M. Ruiz-Altaba and G. Sierra,
	{\it 	
		Quantum Groups in Two-Dimensional Physics (Cambridge Monographs on Mathematical Physics)}
	Cambridge Monographs on Mathematical Physics, Université de Genève,
	Cambridge University Press
	(2010),
	Book ISBN:9780511628825.
	

		\bibitem{Teschner:2002vx}
	J.~Teschner,
	{\it Quantum Liouville theory versus quantized Teichmuller spaces,}
	Fortsch. Phys. \textbf{51} (2003), 865-872, [arXiv:hep-th/0212243].
	%
	\bibitem{Teschner:2003em}
	J.~Teschner, {\it {On the relation between quantum 
			Liouville theory and the quantized Teichm\"uller spaces}}, Int.\ J.\ Mod.\ Phys.\ A. \textbf{19S2} (2004), 459, [arXiv:hep-th/0303149].
	%
	\bibitem{T03}
	J.~Teschner, {\it {A Lecture on the Liouville vertex operators}},  {
		Int. J. Mod. Phys.} \textbf{A19S2} (2004), 436-458, [arXiv:hep-th/0303150].
	%
	\bibitem{Teschner:2005bz}
	J.~Teschner, {\it {An Analog of a modular functor from quantized 
			Teichm\"uller theory}},  in:
	“Handbook of Teichmüller theory”, (A. Papadopoulos, ed.) Volume I, EMS Publishing
	House, Zürich 2007, 685–760, [arXiv:math/0510174].
	%
	\bibitem{teschner} 
	J.~Teschner, {\it {From Liouville theory to the quantum geometry of Riemann surfaces}},   
	Part of Mathematical physics. Proceedings, 14th International Congress, ICMP 2003, Lisbon, Portugal, [arXiv:hep-th/0308031].
	%
	\bibitem{Faddeev:2000if}
	L.~D.~Faddeev, R.~M.~Kashaev and A.~Y.~Volkov,
	{\it Strongly coupled quantum discrete Liouville theory. 1. Algebraic approach and duality},
	Commun. Math. Phys. \textbf{219} (2001), 199-219, [arXiv:hep-th/0006156].
	%
	\bibitem{Faddeev:2008xy}
	L.~D.~Faddeev and A.~Y.~Volkov,
	{\it Discrete evolution for the zero-modes of the Quantum Liouville Model},
	J. Phys. A \textbf{41} (2008), 194008, [arXiv:0803.0230].
	%
	\bibitem{Ka0}
	R.~M. {Kashaev}, {\it 
		Discrete Liouville equation and Teichmüller theory}, [arXiv:0810.4352].
	%
	\bibitem{Ka1}
	R.~M. {Kashaev}, {\it {On the Spectrum of Dehn Twists in Quantum
			Teichm{\"u}ller Theory}},  Physics and Combinatorics, (2001), 63-81, [arXiv:math/0008148].
	%
	\bibitem{Ka2}
	R.~M.~Kashaev, {\it{The quantum dilogarithm and Dehn twists in quantum \TM theory}}, ``Integrable structures of exactly solvable two-dimensional models of quantum field theory'', (Kiev, 2000) 211-221, NATO Sci. Ser. II Math. Phys. Chem., 35, Kluwer Acad. Publ., Dordrecht, 2001.
	%
	\bibitem{Ka4}
	R.~M.~Kashaev, {\it Quantization of \TM spaces and the quantum dilogarithm}, Lett.\ Math.\ Phys.\  \textbf{43} (1998), 105, [arXiv:q-alg/9705021].
	%
	\bibitem{Chekhov:1999tn}
	L.~Chekhov and V.~V.~Fock, {\it Quantum Teichm\"uller space},
	Theor.\ Math.\ Phys.\  \textbf{120} (1999), 1245, [arXiv:math/9908165].
	
	
	
	\bibitem{Peng:2022}
B. Peng, S. Gulania, Y. Alexeev, and N. Govind, \textit{Quantum time dynamics employing the Yang-Baxter equation for circuit compression}, Phys. Rev. A \textbf{106}  (2022) no.1, 012412,
%
	
	
		\bibitem{Kauffman:2004}
	L. H. Kauffman and S. J. Lomonaco, \textit{Braiding operators are universal quantum gates}, New J. Phys. \textbf{6}, 134 (2004),
	[arXiv:quant-ph/0401090].
	%
	
	
	\bibitem{Kitaev:2006}
	A. Kitaev, \textit{Anyons in an exactly solved model and beyond}, Ann. Phys. (N. Y.) \textbf{321}, 2 (2006).
	[arXiv:cond-mat/0506438].
	
	
	

	
	
	
	\bibitem{Kasirajan:2021}
	V. Kasirajan, \textit{Fundamentals of Quantum Computing} (Springer International Publishing, 2021).
	

\bibitem{Aravanis:2022shq}
C.~Aravanis, G.~Korpas and J.~Marecek,
\textit{Transpiling Quantum Circuits using the Pentagon Equation},
[arXiv:quant-ph/2209.14356].



\bibitem{Nill:1996dv}
F.~Nill,
{\it On the structure of monodromy algebras and Drinfeld doubles},
Rev. Math. Phys. \textbf{9} (1997), 371-395,
[arXiv:q-alg/9609020].

	\bibitem{Kulish:1991}
Kulish, P.P. 
{\it Quantum superalgebra $osp(1|2)$.}
J Math Sci \textbf{54}, (1991) 923–930.


\bibitem{Kulish:1989sv}
P.~P.~Kulish and N.~Y.~Reshetikhin,
{\it Universal R matrix of the quantum superalgebra $osp(2 | 1)$},
Lett. Math. Phys. \textbf{18} (1989), 143-149.



\bibitem{Schomerus:2005bf}
V.~Schomerus and H.~Saleur,
{\it The $GL(1|1)$ WZW model: From supergeometry to logarithmic CFT},
Nucl. Phys. B \textbf{734} (2006), 221-245,
[arXiv:hep-th/0510032].



	\bibitem{PT1}
	B.~Ponsot and J.~Teschner,
	{\it Liouville bootstrap via harmonic analysis on a noncompact quantum group},
	[arXiv:hep-th/9911110].
	%
	
	

	\bibitem{PT2} 
	B.~Ponsot and J.~Teschner,
	{\it {Clebsch-Gordan and Racah-Wigner coefficients for a continuous series of representations of $U_q(sl(2,R))$}},	Commun.\ Math.\ Phys.\  \textbf{ 224},  (2001), 613, [arXiv:math/0007097].
%
\bibitem{Hadasz:2013bwa}
L.~Hadasz, M.~Pawelkiewicz and V.~Schomerus,
{\it Self-dual Continuous Series of Representations for $U_q(sl(2))$ and $U_q(osp(1|2))$},
JHEP \textbf{1410} (2014), 91, [arXiv:1305.4596].
%
\bibitem{Pawelkiewicz:2013wga}
M.~Pawelkiewicz, V.~Schomerus and P.~Suchanek,
{\it The universal Racah-Wigner symbol for $U_q(osp(1|2))$}
JHEP \textbf{1404} (2014), 079, [arXiv:1307.6866].
%
	\bibitem{APT}
N.~Aghaei, M.~Pawelkiewicz and J.~Teschner,
{\it Quantisation of super Teichmüller theory},
Commun.\ Math.\ Phys.\  \textbf{353} (2017) no. 2,  597, [arXiv:1512.02617].
%












\bibitem{Andersen:2018pnw}
J.~E.~Andersen and R.~Kashaev,
\textit{The Teichm{\"u}ller TQFT}, Part of Proceedings, International Congress of Mathematicians (ICM 2018): Rio de Janeiro, Brazil, August 1-9, 2018, [arXiv:1811.06853].



\end{thebibliography}
\end{document}